\documentclass[author-year, review]{elsarticle}
\usepackage[margin=3cm]{geometry}
\journal{European Journal of Operational Research}
\usepackage{lineno}
\usepackage{hyperref}

\usepackage{setspace}
\usepackage{color,graphicx,shortvrb}
\usepackage{multirow}
\usepackage{verbatim}
\usepackage{amsmath,amsfonts,amssymb}
\usepackage{rotating}
\usepackage{lscape}
\usepackage{pdflscape}
\usepackage{booktabs}
\usepackage[titletoc,title]{appendix}
\usepackage{hhline}
\usepackage{eurosym}
\usepackage{tikz}
\usepackage{tikz-qtree}
\usepackage{float}
\usepackage{graphics}
\usepackage{subcaption}
\usepackage{changepage}
\usepackage{soul}
\soulregister\cite7
\soulregister\ref7
\soulregister\eqref7
\soulregister\pageref7
\usepackage{cancel}
\usepackage{amsmath,amsfonts,amssymb}
\usepackage{enumitem}
\usetikzlibrary{arrows.meta}
\usepackage{soul}
\usepackage{booktabs}

\usepackage{algorithm,algorithmic}
\makeatletter
\renewcommand{\ALG@name}{Pseudocode}
\makeatother

\usepackage{tikz-cd}
\usepackage{cases}
\usetikzlibrary{matrix,arrows,decorations.pathmorphing}
\usepackage{forest}
\usetikzlibrary{shadows,arrows.meta}

\usepackage{numcompress}\biboptions{authoryear}

\makeatletter
\def\ps@pprintTitle{%
 \let\@oddhead\@empty
 \let\@evenhead\@empty
 \def\@oddfoot{}%
 \let\@evenfoot\@oddfoot}
\makeatother

\usepackage{algorithm,algorithmic}
\makeatletter
\renewcommand{\ALG@name}{Pseudocode}
\makeatother

\newcommand{\be}{\begin{equation}}
\newcommand{\ee}{\end{equation}}
\newcommand{\beenum}{\begin{enumerate}}
\newcommand{\eenum}{\end{enumerate}}
\newcommand{\bi}{\begin{itemize}}
\newcommand{\ei}{\end{itemize}}

\definecolor{darkorange}{rgb}{1.0, 0.55, 0.0}

\newtheorem{theorem}{Theorem}[section]

\usepackage{amsmath}
\tikzset{
dot/.style={circle,draw,inner sep=1.2,fill=black},
}

\newtheorem{proposition}{Proposition}

\begin{document}
 \begin{frontmatter}
	\title{A rolling horizon heuristic approach for a multi-stage stochastic waste collection problem}
		
		\author[SpinelliMaggioni_address]{Andrea Spinelli}
		
		\author[SpinelliMaggioni_address]{Francesca Maggioni\corref{mycorrespondingauthor}}
		\ead{francesca.maggioni@unibg.it}
		\cortext[mycorrespondingauthor]{Corresponding author}
		
		\author[RamosPovoa_address]{T\^ania Rodrigues Pereira Ramos}
		
		\author[RamosPovoa_address]{Ana Paula Barbosa-P\'ovoa}
		
		\author[Vigo_address]{Daniele Vigo}
		
		\address[SpinelliMaggioni_address]{{Department of Management, Information and Production Engineering, University of Bergamo, Viale G. Marconi 5, Dalmine 24044, Italy}}
		\address[RamosPovoa_address]{Centre for Management Studies (CEGIST), Instituto Superior T\'ecnico, University of Lisbon, Av. Rovisco Pais 1, Lisbon 1049-001, Portugal}
		\address[Vigo_address]{Department of Electrical, Electronic, and Information Engineering ``G. Marconi'', University of Bologna, Viale del Risorgimento 2, Bologna 40136, Italy}

		\begin{abstract}
In this paper we present a multi-stage stochastic optimization model to solve an inventory routing problem for recyclable waste collection. The objective is the maximization of the total expected profit of the waste collection company. The decisions are related to the selection of the bins to be visited and the corresponding routing plan in a predefined time horizon. Stochasticity in waste accumulation is modeled through scenario trees generated via conditional density estimation and dynamic stochastic approximation techniques. The proposed formulation is solved through a rolling horizon approach, providing a worst-case analysis on its performance. Extensive computational experiments are carried out on small- and large-sized instances based on real data provided by a large Portuguese waste collection company. The impact of stochasticity on waste generation is examined through stochastic measures, and the performance of the rolling horizon approach is evaluated. Some managerial insights on different configurations of the instances are finally discussed.
		\end{abstract}
		
		
	\end{frontmatter}
	
	
	\section{Introduction} \label{sec_introduction}

In recent years, the importance of sustainable waste management processes has been recognized worldwide (see, for example, the new Circular Economy Action Plan, \cite{Eu2020}). These practices involve long-, medium-, and short-term planning decisions (see \cite{GhiLagManMusVig2014}), and combine different aspects. Among them, the efficiency of the waste collection system is a crucial problem that needs to be addressed (see \cite{BinBloRamPovWonVan2016}). Traditionally, the collection of recyclable household waste (paper, plastic, metal, and glass) is based on static and pre-defined routes, executed on a regular basis (see \cite{MorJorAguPovAntRam2022}). Since this operation does not take into account the actual bins' filling level, high rates of resources' inefficiencies may occur, due to too early collection of not filled bins, or to poor service level because of too late collections. To tackle this problem, some companies have installed sensors inside the waste bins with the aim of tracking the filling rates (see \cite{GutJenHenRia2015,MarBieCorMakVar2020,JorAntRamPov2022}). Such a procedure allows formulating a statistical model of the real amount of waste for each location (see \cite{LopRam2023}). However, a considerable investment is required by service providers to equip all the waste bins of an area with this technology (see \cite{HesDraDoeVig2023}). Therefore, the design of novel Operational Research (OR) techniques is still an important ongoing research topic.

Vehicle Routing Problem (VRP, see \cite{TotVig2002}) is the classical OR modeling approach to tackle waste collection problems (see \cite{HesDraDoeVig2023}). Among all the possible VRP variants, in this paper we focus our attention on the Inventory Routing Problem (IRP, see \cite{CoeCorLap2014b}). According to \cite{MesSchPer2014}, IRPs are particularly effective in waste collection problems since they are characterized by a large number of waste containers, variability in the accumulation rate and a medium to long planning horizon. 

In the great majority of the literature on waste collection models, uncertainty factors are ignored and all the parameters are assumed to be known when making decisions (see \cite{SarGha2023}). Nevertheless, this assumption may not be true in all cases, since uncertainties occur on the traveling time as well as on the accumulation rate of waste (see \cite{MorJorAguPovAntRam2022}). In such a complex framework, Stochastic Optimization techniques (see \cite{BirLou2011}) may help service providers to implement cost-effective decision plans.

Motivated by the uncertain and dynamic nature of the waste accumulation, in this paper we formulate a multi-stage mixed-integer linear stochastic optimization model to solve an Inventory Routing Problem for recyclable waste collection. The waste operator company is required to make decisions at tactical level in a mid-term time horizon. Inventory decisions are integrated in the routing scheme and the planning is based on the amount of waste inside the containers. The aim of the planning is to maximize the profit, given by the difference between the revenues from selling the collected waste and the transportation costs.

This type of problems are among the most challenging in the literature, combining stochasticity and discrete decisions. Exact solution methods are in general based on branch and bound type algorithms or branch and price methods. Since the size of stochastic optimization model grows exponentially with respect to the number of stages and of scenarios, heuristic algorithms are needed. On this purpose, we adopt the \emph{rolling horizon} approach (see \cite{ChaHsuSet2002}), a classical heuristic applied to tackle multi-stage stochastic problems. According to this technique, the model is decomposed into a sequence of subproblems defined over a reduced time horizon. The model is solved starting from the first time period and the value of the first-stage variables is captured. The procedure is then repeated starting from the second stage until the end of the time horizon.

The proposed stochastic formulation is tested on instances of different sizes, based on real data, and the results are validated by means of classical stochastic measures.

The main contributions of this paper can be summarized as follows:

\begin{itemize}
	\item to develop a multi-stage stochastic optimization model for the waste collection inventory routing problem;
	\item to apply the rolling horizon approach to solve the model and analyze its worst-case performance;
	\item to provide numerical experiments with the aim of:
	\begin{enumerate}
		\item[(1)] validating the model in terms of \emph{in-sample stability} (see \cite{KauWal2007});
		\item[(2)] measuring the impact of uncertainty and the quality of the deterministic solution in a stochastic setting;
		\item[(3)] evaluating the performance of the rolling horizon approach in terms of optimal objective function value and reduction of CPU time;
		\item[(4)] testing the effectiveness of the proposed methodology on a real case study.
	\end{enumerate}
\end{itemize}

The remainder of the paper is organized as follows. Section \ref{sec_literature_review} reviews the existing literature on the problem. In Section \ref{sec_problem_description}, the waste collection problem is described and a multi-stage stochastic programming model is formulated. Section \ref{sec_rolling_horizon} describes the rolling horizon approach and provides a worst-case analysis on its performance. In Section \ref{sec_computational_results}, the computational results are shown and the managerial insights are discussed. Finally, Section \ref{sec_conclusions} concludes the paper.

\section{Literature review} \label{sec_literature_review}
Waste collection problems are mostly modeled in the literature as Vehicle Routing Problems, where a predefined set of collection sites is considered and vehicle routes are planned accordingly. The objective is either to minimize transportation cost, total distance, time travelled or to maximize profits, revenue, amount of waste collected (see \cite{HesDraDoeVig2023}). In recent years, such a kind of problems have been widely studied and extended, in order to include different features. To name a few, in \cite{AngSpe2002}, a Periodic Vehicle Routing Problem (PVRP) is designed such that the visiting schedules on a given time horizon are associated with each collection site; \cite{FacPerZan2011} proposes a Capacitated Vehicle Routing Problem (CVRP) where garbage trucks have limited carrying capacity; in \cite{HemDoeHarVig2014}, the problem of designing an integrated system combining PVRP with bins allocation is considered; in \cite{RamMorPov2018}, a CVRP with Profits (CVRPP) is developed such that revenues come from selling the collected waste to a recycling company. Recently, in \cite{TraNguLePhu2024} a novel Location-Assignment-Routing Problem (LARP) is proposed with the aim of optimally locating waste storages, and determining the optimal set of routes in rural developing countries. The reader is referred to \cite{HesDraDoeVig2023} for an updated survey on waste collection routing problems.

Inventory Routing Problem is a modeling extension of VRP because it integrates inventory management and vehicle routing decisions over a medium- or long-term planning period. In the classical IRP approach, three different types of decisions have to be made: when to restock the customers' inventories according to their demand, how much product to deliver, and how to combine customers into vehicle routes. In the special case of waste collection, the flows are reversed because the aim of visiting is collecting rather than delivering (see \cite{MesSchPer2014}). A recent review on routing problems in Reverse Logistics (RL) can be found in \cite{SarGha2023}. Typically, the inventory decisions modelled in IRPs are at customers' locations. In waste collection operations this choice implies that the amount of waste at the collection sites is known and periodic schedules are thus planned (see \cite{HesDraDoeVig2023}). According to \cite{MalSow2018}, IRP models in RL are mostly motivated by real case studies, providing solutions for the management of specific waste types: municipal solid waste (\cite{ElbWoh2016,MorRamPov2019,MarBieCorMakVar2020}), vegetable oil (\cite{AksKaySalAkc2012,AksKaySalTun2014,CarGonTreGar2019,CarMel2021}), infectious medical waste (\cite{NolAbsFeil2014}), end-of-life vehicles (\cite{KriBlaKriFle2008}).

In the traditional VRP approach, all the parameters of the model are considered deterministic (see \cite{GenLapSeg1996}). Nevertheless, in the waste collection system a high degree of uncertainty may affect waste production, disrupting the reliability of the solution made by service providers through deterministic approaches. As an attempt to reduce uncertainty, in some areas waste containers are equipped with volumetric sensors that communicate their filling level to the waste manager (see \cite{RahMohSheHajFusCol2023}). Basing on the transmitted real-time information, the collection is thus dynamically planned. In \cite{RamMorPov2018,JorAntRamPov2022,MorRamLopPov2024}, the Smart Waste Collection Routing Problem (SWCRP) is explored, combining the sensors' usage with optimization techniques. \cite{FadGobPerRosTad2018} derives a scheduling of weekly waste collection activities with sensors installed both on underground containers and inside garbage trucks. However, equipping all the collection sites with sensors results in high financial investments. Even when a partial coverage is enough, additional issues arise in order to find the best locations to place such technology (see \cite{LopRam2023}).

Whenever the supplier has access to some information about the probability distribution of customer's demand, the IRP falls within the framework of Stochastic Inventory Routing Problem (SIRP, see \cite{MoiSal2007,BerBosGueLag2013,CoeCorLap2014b,Ala-OrtDoe2023}). Robust Optimization (RO, see \cite{Ben-TalElGNem2009,SolCorLap2012}) and Chance-Constrained Programming (CCP, see \cite{NemSha2006}) are two alternative paradigms recently explored in the literature to cope with uncertainty in IRP models. Examples of CCP models specifically designed for RL operations under demand uncertainty are discussed in \cite{SoyBloHaiVor2018} and \cite{ZhouLiGuZhaoZieZheFan2021}, whereas in \cite{GhoGohFazMam2022} a RO model for plastic recycling with uncertain demand and transportation cost is formulated. As a downside, including uncertainty in vehicle routing models increases dramatically their computational tractability. To this extent, approximate solution techniques (metaheuristic, heuristic or hybrid algorithms) have been devised in the VRP literature, especially when solving large instances. In the following, we limit our attention to optimization techniques adapted to solve waste collection routing problems under uncertainty.

\cite{NuoKytNisBra2006} considers a stochastic PVRP with time windows for the collection of solid waste in Finland. The solution strategy consists in finding the routing scheme through a combination of a hybrid insertion heuristic and a Guided Variable Neighborhood Thresholding (GVNT) metaheuristic. Waste bins are then sorted by urgency according to their filling rate and scheduling operations are organized. In \cite{KuoZulSur2012}, a CVRP with fuzzy demand is solved with an hybrid metaheuristic combining Particle Swarm Optimization (PSO) and Genetic Algorithm (GA). In \cite{MesSchPer2014}, a stochastic and dynamic IRP is proposed, thanks to the implementation of sensors. Waste bins are categorized into three classes, depending on their waste level and a problem-based heuristic is developed to decide which containers to pick up and on which days. A two-stage stochastic optimization model is formulated in \cite{NolAbsFeil2014} to solve a collection problem of infectious medical waste. Uncertainty lies in the bins' filling rate located at pharmacies and an Adaptive Large Neighborhood Search (ALNS) is used to tackle the complexity of the problem. In \cite{ElbWoh2016}, a multi-product multi-period IRP is considered to collect two different commodities (glass and paper). From historical data, high fluctuations have been observed in waste accumulation rates, leading to consider them as stochastic. A classification procedure is used to categorize containers and a short planning period is defined in a Rolling Horizon (RH) framework. The collection plan is then improved through a Variable Neighborhood Search (VNS) algorithm. \cite{GruFikHuaHirCon-Bol2017} considers a multi-depot VRP with stochastic waste levels. The problem is solved through a combination of an oriented randomization of Iterated Local Search (ILS) and Monte Carlo simulation. In \cite{MarBieCorMakVar2020}, historical data and forecasting techniques are used to estimate the expected containers' filling rate over the planning horizon and to derive the distribution of the overflow probability. Then, an ALNS procedure with search guiding principle based on Simulated Annealing (SA) is developed, in combination with a RH approach. \cite{AliBarPish2021} formulates a bi-objective vehicle routing model with hard time windows for the collection of municipal waste. The amount of generated waste is supposed to be uncertain and modeled through a fuzzy approach. Small instances solutions are generated with the AUGMented $\varepsilon$-CONstraint (AUGMECON) method, while in larger size cases the metaheuristic of Non-dominated Sorting Genetic Algorithm II (NSGA-II) is used. In \cite{ZhouLiGuZhaoZieZheFan2021}, a bi-objective VRP with simultaneous pickup and delivery is considered. The aim is to minimize the total cost and, at the same time, to maximize the recycling revenue. Uncertainty lies in the pick-up demand and is modelled through a fuzzy CCP approach, with the application of Genetic Algorithm (GA) as solution technique. \cite{GhoGohFazMam2022} formulates a robust IRP to model a sustainable RL supply-chain for polystyrene under demand and costs uncertainty. The quality of the solution is then assessed by means of three different heuristics (GA, SA, and Cross Entropy, CE).

All the approaches discussed so far and related to waste collection routing problems under uncertainty are reported in Table \ref{tab_waste_literature_review}. The classification is inspired by \cite{SarGha2023}.

\begin{table}[h!]
\centering
\resizebox{\textwidth}{!}{
    \begin{tabular}{lllllll}
        \toprule
        {Reference} & {Variant} & {Objective} & {Model type} & {Source} & {Solution} & {Waste type}\\
        & {of VRP} & & & {of uncertainty} & {technique}\\
        \hline
        \\ [-1ex]
        \multirow{2}{*}{{\cite{NuoKytNisBra2006}}} & {Stochastic Periodic VRP}  & \multirow{2}{*}{{Minimize distance}} & \multirow{2}{*}{{MILP}} & {Bins' filling rate} & {Heuristic} & \multirow{2}{*}{{Solid waste}}\\
& {with Time Windows} & & & {Travel times} & {(GVNT)}
\\ [1.5ex]
        \multirow{2}{*}{{\cite{KuoZulSur2012}}} & {Capacitated VRP}  & \multirow{2}{*}{{Minimize distance}} & \multirow{2}{*}{{ILP}} & \multirow{2}{*}{{Bins' filling rate}} & {Heuristic} & \multirow{2}{*}{{Solid waste}}\\
& {with Fuzzy Demand} & & & & {(PSO+GA)}
\\ [1.5ex]
\multirow{2}{*}{{\cite{MesSchPer2014}}} & \multirow{2}{*}{{Stochastic Dynamic}  {IRP}} & \multirow{2}{*}{{Minimize cost}} & \multirow{2}{*}{{MILP}} & {Bins' filling rate} & {Heuristic} & \multirow{2}{*}{{Solid waste}} \\
& & & & {Travel time} & {(Problem-based)}
\\ [1.5ex]
\multirow{2}{*}{{\cite{NolAbsFeil2014}}} & \multirow{2}{*}{{Stochastic IRP}} & \multirow{2}{*}{{Minimize cost}} & \multirow{2}{*}{{MILP}} & \multirow{2}{*}{{Bins' filling rate}} & {Heuristic} & \multirow{2}{*}{{Medical products}} \\
& & & & & {(ALNS)}
\\ [1.5ex]
\multirow{2}{*}{{\cite{ElbWoh2016}}} & {Multi-Product} & \multirow{2}{*}{{Minimize cost}} & \multirow{2}{*}{{MILP}} & \multirow{2}{*}{{Bins' filling rate}} & {Heuristic} & \multirow{2}{*}{{Solid waste}}\\
& {Multi-Period IRP} & & & & {(VNS+RH)}
\\ [1.5ex]
\multirow{2}{*}{{\cite{GruFikHuaHirCon-Bol2017}}} & {Multi-Depot VRP} & \multirow{2}{*}{{Minimize cost}} & \multirow{2}{*}{{MILP}} & \multirow{2}{*}{{Bins' filling rate}} & {Heuristic} & \multirow{2}{*}{{Not specified}}\\
& {with Stochastic Demand} & & & & {(ILS+MCS)}
\\ [1.5ex]
{\cite{SoyBloHaiVor2018}} & {Chance-Constrained} {IRP} & {Minimize cost} & {MILP} & {Customer demand} & {Exact} {solver} & {Food waste}\\
\\
\multirow{2}{*}{{\cite{MarBieCorMakVar2020}}} & \multirow{2}{*}{{Stochastic IRP}} & \multirow{2}{*}{{Minimize cost}} & \multirow{2}{*}{{MINLP}} & \multirow{2}{*}{{Bins' filling rate}} & {Heuristic} & \multirow{2}{*}{{Solid waste}}\\
& & & & & {(ALNS+SA+RH)}
\\ [1.5ex]
\multirow{2}{*}{{\cite{AliBarPish2021}}} & {Capacitated VRP}&{Minimize cost} {and} & \multirow{2}{*}{{MILP}} & \multirow{2}{*}{{Bins' filling rate}} & {Heuristic} & \multirow{2}{*}{{Solid waste}} \\
& {with Hard-Time Windows} & {minimize time} & & & {(AUGMECON+NSGA-II)}
\\ [1.5ex]
\multirow{3}{*}{{\cite{ZhouLiGuZhaoZieZheFan2021}}} & {Chance-Constrained} {VRP} &{Minimize cost} {and} & \multirow{3}{*}{{MILP}} & \multirow{3}{*}{{Pick-up levels}} & {Heuristic} & \multirow{3}{*}{{Textile}} \\
& {with Simultaneous} & {maximize profit} & & & {(GA)}\\
& {Pickup and Delivery}
\\ [1.5ex]
\multirow{2}{*}{{\cite{GhoGohFazMam2022}}} & \multirow{2}{*}{{Robust} {IRP}} & \multirow{2}{*}{{Minimize cost}} & \multirow{2}{*}{{MINLP}} & {Demand} & {Heuristic} & \multirow{2}{*}{{Polystyrene}} \\
 & & & & {Costs} & {(GA-SA-CE)}
\\
        \bottomrule    \end{tabular}}
\caption{A selected literature review on the waste collection routing problem under uncertainty.
\\
Abbreviations: ALNS: Adaptive Large Neighborhood Search; AUGMECON: AUGMented $\varepsilon$-CONstraint; CE: Cross Entropy; GA: Genetic Algorithm; GVNT: Guided Variable Neighborhood Thresholding; ILS: Iterated Local Search; MCS: Monte Carlo Simulation; NSGA-II: Non-dominated Sorting Genetic Algorithm II; RH: Rolling Horizon; SA: Simulated Annealing; VNS: Variable Neighborhood Search.} \label{tab_waste_literature_review}
\end{table}

In this work, we test the performance of the rolling horizon approach for the waste collection problem within the paradigm of SIRP. This approach has been extensively used in the literature (see \cite{ChaHsuSet2002} for a classified bibliography and \cite{MagAllBer2014} for the definition of rolling horizon measures). Among its applications to general transportation problems, we mention the works by \cite{maggioni2009stochastic}, \cite{shen2011lagrangian}, and \cite{CavBerMag2022}. In \cite{BerMag2018}, a worst-case analysis of the rolling horizon approach for a stochastic multi-stage fixed charge transportation problem is provided.

\section{Problem description and formulation} \label{sec_problem_description}
A company is responsible for the collection of a recyclable type of waste in a set of $N$ locations (bins or containers) over a time horizon $\mathcal{T}=\{1,\ldots,T\}$. The collection network is represented as a complete directed graph, defined on a set of vertices $\mathcal{I}=\{0,1,\ldots,N\}$ where 0 denotes the depot. Distances $d_{ij}$ are associated with each arc $(i,j)\in\mathcal{I}\times \mathcal{I}$ in the graph. The company needs to determine at stage $t=1$ which waste bins have to be visited and the visiting sequence for all stages $t \in \mathcal{T''}=\{2,\ldots,T\}$. The choice has to be performed with the aim of maximizing the profit over the whole planning horizon $\mathcal{T}$, defined as the difference between the revenues from the selling of the collected waste and the transportation costs. The collected waste is sold at unit price $R$ and the travelling cost per distance unit is fixed to $C$.

Each bin $i\in\mathcal{I'}=\{1,\ldots,N\}$ has a fixed capacity $E_i$ and in the first stage ($t=1$) it is supposed to be filled at $S_{i}^{init}$ percent of its volume. If we assume that the accumulation rate of waste $\{a_i^{(t)}\}_{t=1}^{T}$ of bin $i$ is a random parameter evolving as a discrete-time stochastic process with support $[0,1]$, then the information structure can be described in the form of a scenario tree. At each stage $t\in\mathcal{T}$, there is an ordered set $\mathcal{N}^t=\{1,\ldots,n,\ldots,n^t\}$ of nodes where a specific realization of the uncertain accumulation rate takes place. At the first stage it is associated a unique node $\mathcal{N}^1=\{1\}$, i.e. the root, whereas the final $n^T$ nodes are the leaves of the scenario tree. At stage $t\in\mathcal{T}''$, each node $n\in\mathcal{N}^t$ is connected to a unique node at stage $t-1$, which is called parent (or ancestor) node $pa(n)$. A path through nodes from the root to a leaf is called scenario. At every stage $t\in\mathcal{T}$, each node $n\in\mathcal{N}^t$ has a probability $\pi^{n}$ to occur, and $\sum_{n\in\mathcal{N}^t}\pi^{n}=1$. We denote the accumulation rate for bin $i$ at node $n\in\mathcal{N}^t$, $t\in\mathcal{T}''$ by $a_i^{n}$. For the sake of illustration, we report in Figure \ref{fig_example_scenario_tree} an example of a three-stage scenario tree.

\begin{figure}[ht!]
\begin{minipage}{0.6\textwidth}
\begin{center}
     \centering
\begin{forest}
  for tree={grow = 0, reversed,
    circle, draw, align=center, top color=white, bottom color=white,
        node options={draw},
    edge={},
            parent anchor=east,
        child anchor=west,
    l sep=2cm,
    s sep=0.7cm,
    if level=3{align=left}{},
  }, 
   [1, 
     [1, edge label={node[midway,above=2pt]{\footnotesize$0.7$}}
       [1, edge label={node[midway, above]{\footnotesize$0.4$}}]
       [2, edge label={node[midway, below]{\footnotesize$0.6$}}]
     ]
     [2, edge label={node[midway, below=2pt]{\footnotesize$0.3$}}
      [3, edge label={node[midway, above]{\footnotesize$0.8$}}]
      [4, edge label={node[midway, below]{\footnotesize$0.2$}}]
     ]
   ] 
\end{forest}
\end{center}
\end{minipage} 
\hfill
\begin{minipage}{0.4\textwidth}
	\begin{center}
		\centering
		\resizebox{0.7\textwidth}{!}{
			\begin{tabular}{ccc}\toprule
				{Stage $t$} & {Node $n$} & {$\pi^n$}\\ \hline
				{1} & {1} & {1}\\ \hline
				{2} & {1} & {0.7}\\
				{2} & {2} & {0.3}\\ \hline
				{3} & {1} & {0.28}\\
				{3} & {2} & {0.42}\\
				{3} & {3} & {0.24}\\
				{3} & {4} & {0.06}\\ \bottomrule
		\end{tabular}}
	\end{center}
\end{minipage} 
\caption{Example of a three-stage scenario tree. On the left: the structure of the scenario tree. The numbers on the branches denote the transition probabilities of the random process from each parent node $pa(n)$ to node $n$. On the right: table with nodes at each stage and the corresponding probabilities.}
\label{fig_example_scenario_tree}
\end{figure}

At each stage $t\in\mathcal{T'}=\{1,\ldots,T-1\}$, we define binary decision variables $x_{ij}^t$ and $y_{i}^t$. The former is related to the activation of the arc $(i,j)$ in period $t+1$. Indeed, at each stage the model plans for the next stage, by reflecting what happens in practice for the scheduling of the resources in a waste collection company. If $x_{ij}^t$ is equal to one, then the arc $(i,j)$ will be traversed by a vehicle, with finite capacity $Q$. All the variables $x_{ij}^t$ are defined on the whole graph. Indeed, we assume that, in the collection period, the vehicle starts at the depot, visits the selected bins and returns to the depot to discharge the waste. As far as it concerns the decision variables $y_{i}^t$, if bin $i\in\mathcal{I}'$ needs to be visited in period $t+1$, then variable $y_i^t$ is equal to one at stage $t$. After the realization of the accumulation rate, the amount of waste collected at bin $i$ is denoted by $w_{i}^{n}$, for $n\in\mathcal{N}^t$, $t\in\mathcal{T}'$. In Figure \ref{fig_example_routing} we provide an example of a planning for an horizon of six days.

\begin{figure}[ht!]
\begin{minipage}{0.4\textwidth}
	\begin{center}
		\resizebox{4.4cm}{!}{
			\begin{tikzpicture}
				[node distance={15mm},
				bin/.style = {draw, circle},
				depot/.style = {draw, rectangle}] 
				\node[depot, label = left:\scriptsize{Depot}] (0) {0}; 
				\node[bin] (1) [above right of=0] {1}; 
				\node[bin] (2) [below right of=1] {2}; 
				\node[bin] (3) [below right of=0] {3}; 
				\node[bin] (4) [below left of=0] {4}; 
				\node[bin] (5) [above left of=0] {5}; 
				
				\path[-latex] (0) edge [bend left] node[midway,right]{\scriptsize $x_{01}^1$} (1);
				\path[-latex] (1) edge [bend left] node[midway,left]{\scriptsize $x_{12}^1$} (2);
				\path[-latex] (2) edge [bend left] node[midway,above]{\scriptsize $x_{20}^1$} (0);
				
				\path[-latex] (0) edge [bend right] node[midway,left]{\scriptsize $x_{05}^4$} (5);
				\path[-latex] (5) edge [bend right] node[midway,left]{\scriptsize $x_{54}^4$} (4);
				\path[-latex] (4) edge [bend right] node[midway,below]{\scriptsize $x_{43}^4$} (3);
				\path[-latex] (3) edge [bend left] node[midway,left]{\scriptsize $x_{30}^4$} (0);
		\end{tikzpicture}}
	\end{center}
\end{minipage} 
\hfill
\begin{minipage}{0.6\textwidth}
	\begin{center}
		\centering
		\resizebox{0.9\textwidth}{!}{
			\begin{tabular}{cccc}\toprule
				Stage $t$ & $x_{ij}^t = 1$ & $y_{i}^t = 1$ & Visiting sequence \\ \hline
				1 & $x_{01}^1,x_{12}^1,x_{20}^1$ & $y_{1}^1, y_{2}^1$ \\
				2 & & & $0,1,2,0$\\
				3\\
				4 & $x_{05}^4,x_{54}^4,x_{43}^4, x_{30}^4$ & $y_{3}^4,y_{4}^4,y_{5}^4$ \\
				5 & & & $0,5,4,3,0$\\
				6\\ \bottomrule
		\end{tabular}}
	\end{center}
\end{minipage} 
\caption{Example of a collection plan with 5 bins. On the left: collection routes for day 2 (bins 1, 2) and day 5 (bins 5, 4, 3). On the right: table with active binary decision variables  $x_{ij}^t$ and $y_{i}^t$ and corresponding visiting sequence.}
\label{fig_example_routing}
\end{figure}

At stage $t\in\mathcal{T}'$ and for nodes $n\in\mathcal{N}^t$, additional decision variables are $f_{ij}^{n}$ representing the waste flow shipped through arc $(i,j)$. We assume that the waste flow outgoing depot is null. Finally, for all the time periods, we denote by $u_{i}^{n}$ the accumulated amount of waste at bin $i$. By avoiding partial collection, when bin $i$ is visited, $u_{i}^{n}$ is null.

Moreover, we define the following notation.\\
\emph{\underline{Sets}:}\\
$\mathcal{I}=\{i:i=0,1,\ldots,N\}$: set of $N$ waste bins and the depot, denoted by $0$;\\
$\mathcal{I'}=\{i:i=1,\ldots,N\}$: set of $N$ waste bins (depot excluded);\\
$\mathcal{T}=\{t:t=1,\ldots,T\}$: set of stages;\\
$\mathcal{T'}=\{t:t=1,\ldots,T-1\}$: set of stages (last stage excluded);\\
$\mathcal{T''}=\{t:t=2,\ldots,T\}$: set of stages (first stage excluded);\\
$\mathcal{N}^1=\{n:n=1\}$: root node at stage 1;\\
$\mathcal{N}^t=\{n:n=1,\ldots,n^t\}$: set of ordered nodes of the tree at stage $t\in\mathcal{T}$.\\
\emph{\underline{Deterministic parameters}:}\\
$C$: travelling cost per distance unit;\\
$R$: selling price of a recyclable material;\\
$Q$: vehicle capacity;\\
$B$: waste density;\\
$M$: Big-M number, i.e. a suitable large constant value;\\
$d_{ij}$: distance between $i\in\mathcal{I}$ and $j\in\mathcal{I}$;\\
$S_{i}^{init}$: percentage of waste on the total volume of bin $i\in\mathcal{I'}$ at the first stage;\\
$E_i$: capacity of bin $i\in\mathcal{I'}$;\\
$pa(n)$: parent of node $n\in\mathcal{N}^t$, $t\in\mathcal{T''}$.\\
\emph{\underline{Stochastic parameters}:}\\
$a_i^{n}$: uncertain accumulation rate of bin $i\in\mathcal{I'}$ at node $n\in\mathcal{N}^t$, $t\in\mathcal{T''}$ (percentage on the total volume of the bin);\\
$\pi^{n}$: probability of node $n\in\mathcal{N}^t$, $t\in\mathcal{T}$.\\
\emph{\underline{Decision variables}:}\\
$x_{ij}^t\in\{0,1\}$: binary variable indicating if arc $(i,j)$ is visited at time $t+1$, with $t\in\mathcal{T'}$ and for $i,j\in\mathcal{I}$, $i\neq j$;\\
$y_{i}^t\in\{0,1\}$: binary variable indicating if waste bin $i\in\mathcal{I'}$ is visited at time $t+1$, with $t \in\mathcal{T'}$;\\
$f_{ij}^{n}\in\mathbb{R}^{+}$: non-negative variable representing the waste flow between $i\in\mathcal{I'}$ and $j\in\mathcal{I}$, $i\neq j$, for $n\in\mathcal{N}^t$, $t\in\mathcal{T''}$;\\
$w_i^{n} \in \mathbb{R}^{+}$: non-negative variable representing the amount of waste collected at bin $i\in\mathcal{I'}$, for $n\in\mathcal{N}^t$, $t\in\mathcal{T''}$;\\
$u_i^{n} \in \mathbb{R}^{+}$: non-negative variable representing the amount of waste at bin $i\in\mathcal{I'}$, for $n\in\mathcal{N}^t$, $t\in\mathcal{T}$.

We propose the following stochastic multi-stage mixed-integer linear programming model $\mathcal{M}$:
\begin{eqnarray}	
	\max & &  R\sum_{t\in\mathcal{T''}} \sum_{n\in\mathcal{N}^t} \pi^{n} \sum_{i\in\mathcal{I'}} w_i^{n} -C\sum_{t\in\mathcal{T'}}\sum_{\substack{i,j\in\mathcal{I}\\i\neq j}}d_{ij}x_{ij}^t \label{obj_fun} \\
	\text{s.t. } & & \sum_{\substack{j\in\mathcal{I}\\j\neq i}} f_{ij}^{n} - \sum_{\substack{j\in\mathcal{I'}\\j\neq i}} f_{ji}^{n} = w_{i}^{n} \qquad i\in\mathcal{I'}, n\in\mathcal{N}^t, t\in\mathcal{T''} \label{constraint_01} \\
&	& f_{ij}^{n} \leq (Q-E_jBa_{j}^{n})x_{ij}^{t-1} \qquad i,j \in \mathcal{I}', i\neq j, n\in\mathcal{N}^t, t\in\mathcal{T''} \label{constraint_02} \\
&	& f_{i0}^{n} \leq Qx_{i0}^{t-1} \qquad i\in\mathcal{I}', n\in\mathcal{N}^t, t\in\mathcal{T''} \label{constraint_03} \\
&	& f_{ij}^{n} \leq Q-w_j^{n} \qquad i,j\in\mathcal{I'},i\neq j, n\in\mathcal{N}^t, t\in\mathcal{T''} \label{constraint_04} \\	
 & & f_{ij}^{n} \geq w_i^{n} -M(1-x_{ij}^{t-1})\qquad i\in\mathcal{I'}, j\in\mathcal{I},i\neq j, n\in\mathcal{N}^t, t\in\mathcal{T''} \label{constraint_05} \\
&	& \sum_{\substack{j\in\mathcal{I}\\j\neq i}} x_{ij}^t = y_i^t \qquad i\in\mathcal{I'}, t\in\mathcal{T'} \label{constraint_06}\\
&	& \sum_{\substack{i\in\mathcal{I}\\i\neq j}} x_{ij}^t = y_j^t \qquad j\in\mathcal{I'}, t\in\mathcal{T'} \label{constraint_07}\\
&	& \sum_{i\in\mathcal{I'}}x_{i0}^t = \sum_{j\in\mathcal{I'}}x_{0j}^t \qquad t \in \mathcal{T'} \label{constraint_08}\\
&	& w_i^{n} \leq E_iBy_{i}^{t-1} \qquad i\in\mathcal{I'}, n\in\mathcal{N}^t, t \in \mathcal{T''} \label{constraint_09} \\
&	& u_i^{n} \leq M(1- y_{i}^{t-1}) \qquad i\in\mathcal{I'}, n\in\mathcal{N}^t, t \in \mathcal{T''} \label{constraint_10} \\
&	& u_i^{n} = E_iBS_i^{init} \qquad i\in\mathcal{I'}, n\in\mathcal{N}^1 \label{constraint_11}\\
&	& u_i^{n} = u_i^{pa(n)}+E_iBa_i^{n}-w_i^{n} \qquad i\in\mathcal{I'}, n\in\mathcal{N}^t, t\in\mathcal{T''} \label{constraint_12}\\
&	& u_i^{pa(n)} \leq \big(1 -a_i^{n}\big)E_iB \qquad i\in\mathcal{I'}, n\in\mathcal{N}^t, t\in\mathcal{T''} \label{constraint_13}\\
&	& x_{ij}^t \in\{0,1\} \qquad i,j\in\mathcal{I}, i\neq j, t\in\mathcal{T'} \label{constraint_14}\\
&	& y_{i}^t \in \{0,1\} \qquad i\in\mathcal{I'}, t\in\mathcal{T'} \label{constraint_15}\\
&	& f_{ij}^{n} \geq 0 \qquad i\in\mathcal{I}', j\in\mathcal{I}, i\neq j, n\in\mathcal{N}^t, t\in\mathcal{T''} \label{constraint_16} \\
&	& w _{i}^{n} \geq 0 \qquad i\in\mathcal{I'}, n\in\mathcal{N}^t, t\in\mathcal{T''} \label{constraint_17}\\
&	& u_{i}^{n} \geq 0 \qquad i\in\mathcal{I'}, n\in\mathcal{N}^t, t\in\mathcal{T} \label{constraint_18} 
\end{eqnarray}

The objective function \eqref{obj_fun} is composed by the following terms: (i) the revenues from selling the expected collected waste and (ii) the transportation costs, depending on the routing plan and on the total travelled distance. Constraints \eqref{constraint_01} guarantee the flow balance at each waste bin $i$, for every node $n\in\mathcal{N}^t$ and for every period $t\in\mathcal{T''}$. Constraints \eqref{constraint_02} to \eqref{constraint_04} provide upper bounds on the flow variables $f_{ij}^{n}$, for each node $n\in\mathcal{N}^t$ at stage $t\in\mathcal{T}''$. Specifically, constraints \eqref{constraint_02} guarantee that if bins $i$ and $j$ are not connected, then the waste flow between them is null; otherwise, its sum with the uncertain accumulation amount of waste at $j$ cannot exceed the vehicle capacity. Similarly for constraints \eqref{constraint_03} as far as it concerns the flow between bin $i$ and the depot, once the arc $(i,0)$ is traversed: the vehicle cannot transport to the depot more waste than its capacity. Finally, constraints \eqref{constraint_04} ensure that the sum of the waste flow between bins $i$ and $j$ and the amount of waste collected at bin $j$ cannot exceed the vehicle capacity. Constraints \eqref{constraint_05} provide lower bounds on the flow variable $f^{n}_{ij}$ such that if the vehicle travels from bin $i$ to bin $j$ or from bin $i$ to the depot, with $n\in\mathcal{N}^t$ and $t\in\mathcal{T}''$, all of the accumulated amount of waste at bin $i$ should be collected. Constraints \eqref{constraint_06} and \eqref{constraint_07} link together the decision variables $x_{ij}^t$ and $y_{i}^t$ for each stage $t\in\mathcal{T'}$ and ensure that, if bin $i$ is visited, then there exists exactly one route reaching and one route leaving $i$; on the other hand, no visits at bin $i$ imply no incoming edges to and no outgoing edges from $i$. Constraints \eqref{constraint_08} impose the depot's balance by enforcing that the numbers of incoming and outgoing edges are the same for every period $t\in\mathcal{T'}$. This means that, whether the vehicle performs a route starting from the depot, then it must return to the depot. Constraints \eqref{constraint_09} ensure that the collection amount $w_i^{n}$ at bin $i$ in node $n\in\mathcal{N}^t$, for $t\in\mathcal{T}''$ must be zero, unless the bin is visited. Constraints \eqref{constraint_10} guarantee that the amount of waste $u_i^{n}$ at bin $i$ at node $n\in\mathcal{N}^t$ and stage $t\in\mathcal{T''}$ must be zero if the bin is visited. Constraints \eqref{constraint_11} fix the initial amount of waste $u_i^{n}$ at bin $i$ at the root of the scenario tree. Constraints \eqref{constraint_12} update at every node $n\in\mathcal{N}^t$ and for every period $t\in\mathcal{T''}$ the amount of waste $u_i^{n}$ at bin $i$ by incorporating the uncertain accumulated amount of waste and, potentially, by subtracting the amount of collected waste $w_i^{n}$. Constraints \eqref{constraint_13} impose that no bins are allowed to overflow at each node $n\in\mathcal{N}^t$, $t\in\mathcal{T}''$. Finally, constraints from \eqref{constraint_14} to \eqref{constraint_18} define the decision variables of the problem. We denote by $z^{*}$ the optimal expected profit of model $\mathcal{M}$.

\subsection{A two-commodity flow model}
In model $\mathcal{M}$, the distances between two locations are considered as asymmetric, i.e. in general $d_{ij}\neq d_{ji}$, for $i,j\in\mathcal{I}$. This assumption impacts not only the objective function \eqref{obj_fun}, but also both the constraints \eqref{constraint_01} and \eqref{constraint_03}-\eqref{constraint_05} related to the flow variables $f_{ij}^n$, and the degree constraints \eqref{constraint_06}-\eqref{constraint_08} on variables $x_{ij}^t$ and $y_i^t$. This leads to an increase of the size of the model, due to a considerable number of inequality constraints.

In practical cases, however, distances $d_{ij}$ and $d_{ji}$ on arcs $(i,j)$ and $(j,i)$, respectively, may not be significantly different. For example, in the largest instance of this work the average difference $\lvert d_{ij}-d_{ji}\rvert$ is of hundreds of meters, corresponding to an average percentage difference of 6\%. Therefore, considering a symmetric distance matrix does not result in a considerable worsening of the solution. For this reason, we design an alternative version of model $\mathcal{M}$, denoted by $\mathcal{M}_{sym}$, based on the two-commodity flow formulation proposed in \cite{BalHadMin2004} and applied to a waste collection problem in \cite{RamMorPov2018}. Hence, a copy depot denoted by $N+1$ is introduced and each route is defined according to two paths: one direct path from depot $0$ to depot $N+1$, with variables $f_{ij}^n$ representing the load of the vehicle, and one reverse path, from depot $N+1$ to depot $0$, with variables $f_{ji}^n$ denoting the empty space of the vehicle (see Figure \ref{fig_example_routing_twoflow} for an illustrative example).

\begin{figure}[h!]
\begin{center}
	\resizebox{12cm}{!}
	{
		\begin{tikzpicture}
			[node distance={30mm},
			bin/.style = {draw, circle},
			depot/.style = {draw, rectangle}] 
			\node[depot, label = below:\scriptsize{{Real depot}}] (0) {0}; 
		    \node[bin] (5) at (-4,1) {5}; 
			\node[bin] (1) [above right of=0] {1}; 
			\node[bin] (4) [below of=5] {4}; 
			\node[depot, label = above:\scriptsize{{{Copy depot}}}] (6) at (1,0) {{6}}; 
		    \node[bin] (3) at (4,-2) {3}; 
		    			\node[bin] (2) [above right of=3] {2}; 
			
			\path[-latex] (0) edge [bend right=40] node[midway, above=2pt]{\scriptsize $f_{05}^{n}=0$} (5);
			\path[-latex] (5) edge [bend right] node[midway, left=2pt ]{\scriptsize $f_{54}^{n}=3$} (4);
			\path[-latex] (4) edge [bend right=15] node[midway, below=2pt ]{\scriptsize $f_{43}^{n}=4$} (3);
			\path[-latex] (3) edge [bend right=40] node[midway, right=2pt]{\scriptsize {$f_{36}^{n}=7$}} (6);

			\path[-latex, dashed] (5) edge [bend right=16] node[midway, above=2pt]{\scriptsize $f_{50}^{n}=7$} (0);
			\path[-latex, dashed] (4) edge [bend right] node[midway, right=2pt]{\scriptsize $f_{45}^{n}=4$} (5);
			\path[-latex, dashed] (3) edge [bend right=10] node[midway, below]{\scriptsize $f_{34}^{n}=3$} (4);
			\path[-latex, dashed] (6) edge [bend right=20] node[midway,  above right=-0.5pt]{\scriptsize {$f_{63}^{n}=0$}} (3);
	\end{tikzpicture}}
\end{center}
\caption{Representation of the two-commodity flow formulation on the same network of Figure \ref{fig_example_routing}. A copy depot (vertex 6) is introduced, and the truck capacity $Q$ is set to 7. The solid lines represent the actual visiting sequence, starting from the real depot, with corresponding waste flows $f^n_{ij}$. The dashed lines are associated with the reverse flows $f^n_{ji}$, related to the empty space in the vehicle. Note that $f^n_{ij}+f^n_{ji}=Q$.}
\label{fig_example_routing_twoflow}
\end{figure}

Each edge is therefore counted twice and the objective function \eqref{obj_fun} needs to be updated as
\begin{linenomath}
\begin{equation*}
	\max \quad R\sum_{t\in\mathcal{T''}} \sum_{n\in\mathcal{N}^t} \pi^{n} \sum_{i\in\mathcal{I'}} w_i^{n} -{\color{black}{\frac{C}{2}}}\sum_{t\in\mathcal{T'}}\sum_{\substack{i,j\in\mathcal{I}\\i\neq j}}d_{ij}x_{ij}^t.
\end{equation*}
\end{linenomath}

Constraints \eqref{constraint_01} are replaced by
\begin{linenomath}
\begin{equation*}
	\sum_{\substack{j\in\mathcal{I}\\j\neq i}} (f_{ij}^{n} - f_{ji}^{n}) = 2w_{i}^{n}\qquad i\in\mathcal{I'}, n\in\mathcal{N}^t, t\in\mathcal{T''},
\end{equation*}
\end{linenomath}

since the two-commodity flow formulation considers two flows passing through each node $i$. In addition, constraints \eqref{constraint_03}-\eqref{constraint_05} are substituted by
\begin{linenomath}
\begin{equation}\label{constraint_symm_depot}
	\sum_{i\in\mathcal{I'}} f_{iN+1}^{n} = \sum_{i\in\mathcal{I'}}w_{i}^{n} \qquad n\in\mathcal{N}^t, t\in\mathcal{T''}
\end{equation}
\end{linenomath}

and
\begin{linenomath}
\begin{equation}\label{constraint_symm_flows}
	f_{ij}^{n}+f_{ji}^n = Qx_{ij}^{t-1} \qquad i,j \in \mathcal{I}, i\neq j, n\in\mathcal{N}^t, t\in\mathcal{T''}.
\end{equation}
\end{linenomath}

Constraints \eqref{constraint_symm_depot} ensure that the total inflow of the copy depot corresponds to the total amount of collected waste, whereas constraints \eqref{constraint_symm_flows} impose that, whenever an edge is traversed, the sum of the direct and reverse flows is equal to the capacity of the vehicle. Finally, the degree constraints \eqref{constraint_06}-\eqref{constraint_08} reduce to
\begin{linenomath}
\begin{equation*}
	\sum_{\substack{i\in\mathcal{I}\\i\neq j}} x_{ij}^t = 2y_j^t \qquad j\in\mathcal{I'}, t\in\mathcal{T'}.
\end{equation*}
\end{linenomath}

All the other constraints not mentioned remain unchanged when passing from model $\mathcal{M}$ to $\mathcal{M}_{sym}$. For the sake of completeness, the entire model formulation $\mathcal{M}_{sym}$ is reported in the Appendix.

\subsection{A polynomially solvable case}
The proposed SIRP formulation is clearly NP-hard, since it can be reduced to the well-known NP-hard Travelling Salesman Problem (see \cite{GarJho1979}), whenever the time horizon is $\mathcal{T}=\{1,2\}$, the selling price $R$ is zero, the capacity $Q$ of the vehicle is infinite, and all the waste bins need to be visited at day 2 in order to avoid overflow.

On the other hand, the waste collection problem admits a polynomially solvable case if routing decisions are excluded from the problem. This will be addressed in the following proposition.

\begin{proposition} \label{thm_comp_compl_C0}
If $C=0$, then the optimal profit of model $\mathcal{M}$ is
\begin{linenomath}
	\begin{equation}\label{profit_C0}
		z^* = RB\bigg\{\sum_{i\in\mathcal{I}'}E_i\bigg(S_i^{init}+\sum_{t\in\mathcal{T}''}\mathbb{E}\big[a_i^{(t)}\big]\bigg)\bigg\},
	\end{equation}
\end{linenomath}
where $\mathbb{E}\big[a_i^{(t)}\big]$ is the expected accumulation rate of waste at time $t\in \mathcal{T}''$ for bin $i\in\mathcal{I}'$.
\end{proposition}

\noindent{\textit{Proof.}}
We prove the proposition by induction on the time horizon $T$.
\begin{itemize}
	\item (Base case) We consider the case of a two-stage problem ($T=2$). Since $C=0$ and $\mathcal{T}''=\{2\}$, profit \eqref{obj_fun} reduces to
	\begin{linenomath}
		\begin{equation*}
			z=R \sum_{n\in\mathcal{N}^2}\pi^n\sum_{i\in\mathcal{I}'}w_i^n.
		\end{equation*}
	\end{linenomath}
	From constraints \eqref{constraint_11}-\eqref{constraint_12}, it holds that
	\begin{linenomath}
		\begin{equation*}
			w_i^n = E_iBS_i^{init}+E_iBa_i^n-u_i^n \qquad i\in\mathcal{I}',n\in\mathcal{N}^2,
		\end{equation*}
	\end{linenomath}
	which, substituting in the objective function and considering that $\sum_{n\in\mathcal{N}^2}\pi^n=1$ and $\sum_{n\in\mathcal{N}^2}\pi^na_i^n=\mathbb{E}\big[a_i^{(2)}\big]$, gives
	\begin{linenomath}
		\begin{equation*}
			\begin{split}
				z & = R\sum_{n\in\mathcal{N}^2}\pi^n\sum_{i\in\mathcal{I}'}E_iBS_i^{init}+R\sum_{n\in\mathcal{N}^2}\pi^n\sum_{i\in\mathcal{I}'}E_iBa_i^n-R\sum_{n\in\mathcal{N}^2}\pi^n\sum_{i\in\mathcal{I}'}u_i^n\\
				& = RB \sum_{i\in\mathcal{I}'}E_i(S_i^{init}+\mathbb{E}\big[a_i^{(2)}\big])-R\sum_{n\in\mathcal{N}^2}\pi^n\sum_{i\in\mathcal{I}'}u_i^n.
			\end{split}
		\end{equation*}
	\end{linenomath}
Moreover, we note that the objective function $z$ is the difference of two non-negative quantities, where the first one is constant. Thus, 
	\begin{equation*}
		\max z= RB \sum_{i\in\mathcal{I}'}E_i(S_i^{init}+\mathbb{E}\big[a_i^{(2)}\big]) -\min R\sum_{n\in\mathcal{N}^2}\pi^n\sum_{i\in\mathcal{I}'}u_i^n,
		\end{equation*} 
where the minimum of the second term is reached at $u_i^n=0$, for all $n\in\mathcal{N}^2$, $i\in\mathcal{I}'$: the thesis is verified.
	\item (Inductive step) We assume that the thesis holds for a model with time horizon $T-1$, being $T>2$. We need to prove that the thesis is also verified for a model with time horizon $T$. In this case, the objective function can be decomposed as
	\begin{linenomath}
		\begin{equation*}
			R \sum_{t =2}^{T} \sum_{n\in\mathcal{N}^t}\pi^n\sum_{i\in\mathcal{I}'}w_i^n = R \sum_{t =2}^{T-1} \sum_{n\in\mathcal{N}^t}\pi^n\sum_{i\in\mathcal{I}'}w_i^n + R \sum_{n\in\mathcal{N}^{T}}\pi^n\sum_{i\in\mathcal{I}'}w_i^n.
		\end{equation*}
	\end{linenomath}
	Given the induction hypothesis, the optimal profit of the first addendum corresponds to $RB\bigg\{\sum_{i\in\mathcal{I}'}E_i\bigg(S_i^{init}+\sum_{t=2}^{T-1}\mathbb{E}\big[a_i^{(t)}\big]\bigg)\bigg\}$. At stage $T$, from constraints \eqref{constraint_12}, it holds that $w_i^n=E_iBa_i^n$, for all $i\in\mathcal{I}'$, $n\in\mathcal{N}^{T}$, since $u_i^{pa(n)}=0$ for the induction hypothesis and $u_i^n=0$ for the same reasoning of the base case. Consequently, we get
	\begin{linenomath}
		\begin{equation*}
			z^* = RB\bigg\{\sum_{i\in\mathcal{I}'}E_i\bigg(S_i^{init}+\sum_{t=2}^{T-1}\mathbb{E}\big[a_i^{(t)}\big]\bigg)\bigg\} + RB\sum_{i\in \mathcal{I}'} E_i \mathbb{E}\big[a_i^{(T)}\big],
		\end{equation*}
	\end{linenomath}
	which verifies the thesis. \qed
\end{itemize}

We conclude that parameters $R$ and $C$ have different roles: when $R=0$ model $\mathcal{M}$ is NP-hard, whereas if $C=0$ an optimal policy can be computed in $O(T)$ time. For this reason, given the computational complexity of the problem in the general case, heuristic methods are required. To cope with this issue, in the next section we apply the rolling horizon approach to the considered problem.

\section{The rolling horizon approach and its worst-case analysis}\label{sec_rolling_horizon}

One of the most classical heuristic algorithms for multi-stage stochastic programming models is the rolling horizon approach (see \cite{ChaHsuSet2002}). According to this methodology, the multi-stage stochastic problem is decomposed into a sequence of subproblems with a fewer number $W$ of consecutive periods (see \cite{CavBerMag2022}). This leads to a reduced computational effort because at each iteration of the algorithm the number of nodes considered in the scenario tree is lower than the one in the original multi-stage program. However, the quality of the solution may deteriorate since the time horizon is reduced and the solution may be suboptimal (see \cite{BerMag2018}). In the following we present the details of the approach.

First of all, we fix the reduced number $W$ of consecutive period, with $1\leq W < T-1$. In the first iteration of the algorithm, the $(W+1)$-stage stochastic programming model defined on $t=1,\ldots,W+1$ is solved, and the values of the first-stage decision variables $x_{ij}^1$ and $y_i^1$ and the second-stage variables $w_{i}^n$ and $u_{i}^n$, for $n\in\mathcal{N}^2$, are stored. In the second iteration, the value of the inventory levels $u_{i}^n$ for $n\in\mathcal{N}^2$ are fixed as the ones deduced from the first iteration. This is needed to keep track of the evolution of the process and to link the two consecutive time periods. Then, the $(W+1)$-stage stochastic programming model defined on $t=2,\ldots,W+2$ is solved and, as before, the values of the second-stage decision variables $x_{ij}^2$ and $y_i^2$ and the third-stage variables $w_{i}^n$ and $u_{i}^n$, for $n\in\mathcal{N}^3$, are stored. This process is repeated until the last iteration defined on stages $t=T-W,\ldots,T$ is performed. Then, a $W$-stage stochastic programming model defined on $t=T-W+1,\ldots,T$ is solved and the same approach described above is applied. Next, a $(W-1)$-stage stochastic programming model defined on $t=T-W+2,\ldots,T$ is solved and the process is repeated until the last two-stage stochastic programming model defined on $t=T-1,T$.

Once the $T-1$ stochastic programming models have been solved, the variables $x_{ij}^t$, $y_{i}^t$ for all $t\in\mathcal{T}'$ and $w_{i}^n$ for all $n\in\mathcal{N}^t$ and $t\in\mathcal{T}''$ are obtained. The corresponding value of the objective function \eqref{obj_fun} is then computed, leading to $z^{RH,W}$.

Schematically, the algorithm can be represented as in Pseudocode \ref{pseudocode_RH}.

\begin{algorithm}[H]
\begin{algorithmic}[1]
	\REQUIRE $T$, $1\leq W < T-1$
	\STATE $k \gets 1$, $l \gets W+1$
	\STATE $u_i^{+n*} \gets u_i^1$, $n\in\mathcal{N}^t$, $t=k$
	\WHILE{$k\leq T-W$}
	\STATE Solve $(W+1)$-stage SP on $t=k,\ldots,l$ with $u_i^{+n} \gets u_i^{+n*}$, $n\in\mathcal{N}^t$, $t=k$
	\STATE Store $x_{ij}^{t*},y_i^{t*}$, $t=k$ and $u_i^{n*}, w_i^{n*}$, $n\in\mathcal{N}^t$, $t=k+1$
	\STATE $k \gets k+1$, $l \gets l+1$
	\ENDWHILE
	\STATE $j\gets 1$
	\WHILE{$k\leq T-1$}
	\STATE Solve $(W+1-j)$-stage SP on $t=k,\ldots,T$ with $u_i^{+n} \gets u_i^{+n*}$, $n\in\mathcal{N}^t$, $t=k$
	\STATE Store $x_{ij}^{t*},y_i^{t*}$, $t=k$ and $u_i^{n*}, w_i^{n*}$, $n\in\mathcal{N}^t$, $t=k+1$
	\STATE $k \gets k+1$, $j \gets j+1$
	\ENDWHILE
	\STATE Return the corresponding value of the objective function \eqref{obj_fun}.
\end{algorithmic}
\caption{The rolling horizon approach for model $\mathcal{M}$} \label{pseudocode_RH}
\end{algorithm}

We now perform  a worst-case analysis of this approach. The following results hold true.
\begin{theorem} \label{thm_worstcase_RH_C0}
If $C=0$, then
\begin{linenomath}
	\begin{equation*}
		\frac{z^{RH,W}}{z^{*}}=1,
	\end{equation*}
\end{linenomath}
for every choice of $W=1,\ldots,T-2$.
\end{theorem}
\textit{Proof.}
Consider the case of $W=1$, where $T-1$ two-stage stochastic optimization models have to be solved. Since all the subproblems do not share any overlapping period, if we denote by $z^{RH,1}_{t,t+1}$ the objective function value on time period $\{t,t+1\}$, the total profit is
\begin{linenomath}
	\begin{equation*}
		\begin{split}
			z^{RH,1} &= z^{RH,1}_{1,2} + z^{RH,1}_{2,3} + \ldots + z^{RH,1}_{T-1,T}\\
			& = RB\bigg\{\sum_{i\in\mathcal{I}'}E_i\bigg(S_i^{init}+\mathbb{E}\big[a_i^{(2)}\big]+\mathbb{E}\big[a_i^{(3)}\big]+\mathbb{E}\big[a_i^{(T)}\big]\bigg)\bigg\}.
		\end{split}
	\end{equation*}
\end{linenomath}
The previous expression coincides with \eqref{profit_C0}, and so the thesis is verified.
 
When considering a value $W>1$, only the collecting variables $w_{i}^n$ at the second stage of each subproblem are stored, meaning that exclusively the accumulation rates at that stage are considered in the optimal solution. This implies the thesis in a similar fashion as $W=1$. \qed

On the other hand, when $R=0$, the following result on the performance of the rolling horizon approach with $W=1$ holds.

\begin{figure}[h!]
	\centering
	\includegraphics[width=0.9\textwidth]{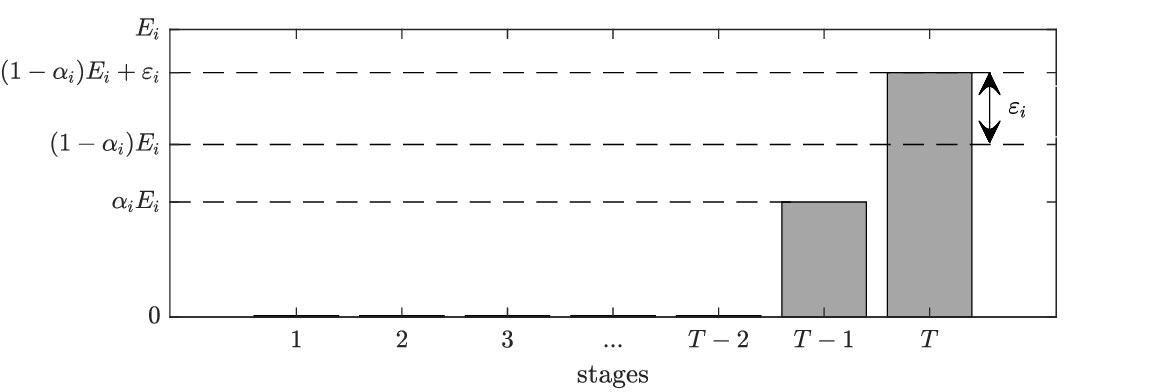}
	\caption{Graph of the accumulation rate \eqref{acc_rate_worstcase}.}
	\label{fig_acc_rate_worst_case}
\end{figure}

\begin{theorem} \label{thm_worstcase_RH_W1}
There exists a class of instances such that $z^{RH,1}=-\infty$, even if model $\mathcal{M}$ is feasible.
\end{theorem}
\textit{Proof.}
Consider the following class of instances: initial amount of waste $S_{i}^{init}=0$ for all $i\in\mathcal{I}'$; vehicle capacity $Q>\sum_{i\in\mathcal{I}'}E_i$; selling price $R=0$.\\
For all $i\in\mathcal{I}'$, let $\alpha_i\in (0;1)$, $\varepsilon_i\in(0;\alpha_i E_i]$, and the accumulation rate $a_i^n$ be such that
\begin{linenomath}
	\begin{equation}
		a_i^n = 
		\begin{cases}
			0 & \text{if } n\in\mathcal{N}^t, \text{ } t\in\mathcal{T}'\cup\{T-2\}\\
			\alpha_iE_i & \text{if } n\in\mathcal{N}^{T-1}\\
			(1-\alpha_i)E_i+\varepsilon_i & \text{if } n\in\mathcal{N}^{T}.
		\end{cases} \label{acc_rate_worstcase}
	\end{equation} 
\end{linenomath}
 The graph of $a_i^n$ is depicted in Figure \ref{fig_acc_rate_worst_case}.

We apply the rolling horizon approach with $W=1$. This means that $T-1$ two-stage stochastic programming models have to be solved. In the first $T-3$ programs, all the decision variables are zero, since there is no waste to collect. Similarly for the model defined on stages $t=T-2,T-1$, with the exception of $u_i^{T-1}$ which is equal to $\alpha_iE_i$, for all $i\in \mathcal{I}'$. However, the last optimization model defined on stages $t=T-1,T$ is infeasible because each waste bin $i\in\mathcal{I}'$ incurs into overflowing, due to the violation of constraints \eqref{constraint_13}. This implies that $z^{RH,1}=-\infty$.

On the other hand, the optimal profit $z^*$ of the multi-stage stochastic optimization model is equal to $-C\sum_{\substack{i,j\in\mathcal{I}\\i\neq j}}d_{ij}x_{ij}^{T-1}$, deriving from a collection on day $T$. \qed

Making the appropriate changes, a similar performance of the rolling horizon approach with $W=1$ also holds within model $\mathcal{M}_{sym}$.

\section{Computational results} \label{sec_computational_results}
In this section, we first describe the instances on which we perform the numerical simulations (see Subsection \ref{sec_data_analysis}). Subsection \ref{sec_analysis_of_the_solution} compares the solutions of models $\mathcal{M}$ and $\mathcal{M}_{sym}$. In Subsection \ref{sec_stochastic_measure}, the validation of model $\mathcal{M}$ with standard stochastic measures is provided, and the quality of the expected value solution is discussed. In Subsection \ref{sec_performance_RH}, the performance of the rolling horizon approach is assessed. Then, the results on a large case study are presented in Subsection \ref{sec_real_case_study}. Finally, managerial insights are provided in Subsection \ref{sec_managerial_insights}.

All computational experiments are obtained using GAMS 38.3.0 and solver Gurobi 9.5 on an Intel(R) Core(TM) i5-8500 64-bit machine, with 8 GB RAM and 3.00 GigaHertz processor. Unless otherwise specified, a runtime limit of 24h is imposed.

\subsection{Data analysis} \label{sec_data_analysis}
The data considered in this study are inspired by a real case problem provided by the industrial partner \textit{ERSUC - Res\'iduos S\'olidos do Centro, S.A.}, one of the main waste management companies in Portugal. The company operates in the Central Region of Portugal and owns a homogeneous fleet of vehicles based at two different depots, one near the city of \textit{Aveiro} and the other close to \textit{Coimbra}. The recyclables collection is performed independently for each type of waste (glass, paper/cardboard and plastic/metal).

The case study described in the following focuses on the collection of plastic/metal waste, related to packaging materials and around the suburban municipality of \textit{Condeixa-a-Nova} in the district of \textit{Coimbra}. Real data have been provided by \textit{ERSUC} on the filling rate of 121 waste bins between April and July 2019 (15 weeks). The data are gathered by the garbage collector only on the collection days (20 days in total). The working days of the company include all the days of a week, except Sunday. Therefore, we set the time horizon in the model as $\mathcal{T}=\{1,\ldots,6\}$.

We perform simulations on both small and large instances. As far as it concerns the small cases, we generate a set of thirty instances with a reduced number of bins, randomly drawn from the entire dataset of 121 bins. For simplicity, we denote each small instance on the basis of the coding scheme ``\emph{inst\_$draw$\_$numbins$}'', where $draw$ is an integer between 1 and 10 associated with the random draw, and $numbins$ is the number of selected bins (9, 10 or 11). In addition, we consider a large instance composed by 50 bins to simulate a real case study of waste collection, since fifty is the average number of bins in a collection route of the industrial partner.

The deterministic parameters of the model are shown in Table \ref{tab_deterministic_parameters}. 

\begin{table}[h!]
\centering
\resizebox{\textwidth}{!}{
	\begin{tabular}{lll}\toprule
		Parameter  & Value & Source \\ \hline
		$C$ & 1 \euro/km & \textit{ERSUC}\\
		$R$ & 0.30 \euro/kg & \textit{Sociedade Ponto Verde}\\
		$Q$ & 2000 kg & \textit{ERSUC}\\
		$B$ & 30 kg/m$^3$ & \textit{ERSUC}\\
		$M$ & $10^5$ & -\\
		$d_{ij}$, $i,j\in\mathcal{I}$ & Actual road distance between $i$ and $j$ & \textit{ERSUC} and \textit{OpenRouteService}\\
		$S_i^{init}$, $i\in\mathcal{I}'$ & Initial percentage of waste on the total volume of bin $i$ & \textit{ERSUC}\\
		$E_i$, $i\in\mathcal{I}'$ & 2.5 m$^3$ & \textit{ERSUC}\\
		\bottomrule
\end{tabular}}
\caption{Parameters values and sources.} \label{tab_deterministic_parameters}
\end{table}

As in \cite{RamMorPov2018}, transportation cost $C$ includes fuel consumption, maintenance of the vehicle and drivers' wages. The revenue parameter $R$ is derived as follows: for each ton of packaging collected and sorted, the \textit{Sociedade Ponto Verde} (the packaging waste regulator in Portugal) pays 545 \euro/ton to the waste collection company; since only the collection activity is being considered in this work, which corresponds to approximately the $55\%$ of the total cost, the selling price $R$ is adjusted to 0.30 \euro/kg.

We discuss now how we construct the random process $\{a_i^{(t)}\}_{t=1}^{T}$ of the daily accumulation rate, based on observations provided by the industrial partner. For each bin $i\in\mathcal{I}'$, we denote by $\{p_{i}^{(t)}\}_{t=1}^{20}$ the historical data of the filling rate on collection days. For all $i\in\mathcal{I'}$, we set $S_{i}^{init}$ equal to $p_i^{(1)}$. We assume that, if $t_{1}$ and $t_{2}$ are two consecutive collection days, the increase (or decrease) of the filling rate of waste between $t_{1}$ and $t_{2}$ is constant. Once the waste collector visits bin $i$ at time $t_{1}$, then she/he empties it, i.e. $p_{i}^{(t_{1})}=0$. Thus, the daily accumulation rate of waste in bin $i$ can be calculated as
\begin{linenomath}
\begin{equation*}
	a_{i}^{(t)}=\frac{p_{i}^{(t_{2})}-p_{i}^{(t_{1})}}{t_{2}-t_{1}} = \frac{p_{i}^{(t_{2})}}{t_{2}-t_{1}}, \qquad t=t_{1}+1,\ldots,t_{2}.
\end{equation*}
\end{linenomath}
By following this procedure, for each bin $i\in\mathcal{I}'$, a complete trajectory of the stochastic process $\{a_i^{(t)}\}_{t=1}^{T}$ is obtained on a daily basis.

In the Appendix we report the scenario tree generation procedure we adopt, along with an in-sample stability analysis on the number of scenarios to be considered in the scenario tree. In the remainder of the paper, we show the results obtained on a tree with 32 scenarios and 63 nodes.

\subsection{A comparison of models $\mathcal{M}$ and $\mathcal{M}_{sym}$ solutions} \label{sec_analysis_of_the_solution}
Solving either model $\mathcal{M}$ or model $\mathcal{M}_{sym}$ to optimality on the whole dataset of 121 bins is not possible on our machine, given the high number of variables and constraints (see Table \ref{tab_summary_results_big_instance}). When considering reduced instances composed by 9, 10 or 11 bins, models $\mathcal{M}$ and $\mathcal{M}_{sym}$ provide the same policy in terms of bin selection, visiting schedule, and consequent weight of collected waste (see row 7 in  Table \ref{tab_summary_results_small_instances}). However, due to the assumption of symmetric distances, the travelled distance and the profit are different. On the other hand, on the large instance with 50 bins model $\mathcal{M}_{sym}$ outperforms model $\mathcal{M}$, since the optimality gap is much smaller (see row 11 in Table \ref{tab_summary_results_big_instance}). For this reason, in the following we show results obtained with model $\mathcal{M}$ on small instances, whereas the large instance outcomes rely on model $\mathcal{M}_{sym}$.

\begin{table}[h!]
\centering
\resizebox{\textwidth}{!}{
	\begin{tabular}{l|ll|ll|ll}\toprule
		Number of bins & 9 ($\mathcal{M}$) & 9 ($\mathcal{M}_{sym}$) & 10 ($\mathcal{M}$) & 10 ($\mathcal{M}_{sym}$) & 11 ($\mathcal{M}$) & 11 ($\mathcal{M}_{sym}$) \\ \hline
		Binary variables & 495 & 595 & 600 & 710 & 715 & 835 \\
		Continuous variables & 6705 & 7263 & 8070 & 8690 & 9559 & 10241\\
		Equality constraints & 1220 & 8052 & 1355 & 9546 & 1490 & 11164 \\
		Inequality constraints & 16182 & 6138 & 19840 & 7440 & 23870 & 8866\\ \hline
		Profit (\euro) & 18.16 {(13.27)} & 18.45 {(13.40)} & 34.37 {(13.40)} & 34.67 {(13.62)} & 41.25 {(14.88)} & 41.{5}0 {(14.81)} \\
		Weight of collected waste (kg) & 478.15 {(102.51)} & 478.15 {(102.51)} & 513.83 {(76.13)} & 513.83 {(76.13)} & 612.48 {(44.04)} & 612.48 {(44.04)} \\
		Travelled distance (km) & 125.28 {(22.08)} & 124.99 {(21.96)} & 119.78 {(19.50)} & 119.46 {(19.33)} & 142.50 {(8.66)} & 142.25 {(8.66)}\\
		Ratio weight/distance (kg/km) & 3.80 {(0.32)} & 3.80 {(0.32)} & 4.32 {(0.45)} & 4.33 {(0.46)} &  4.31 {(0.38)} & 4.32 {(0.38)}\\
		CPU time (s) & 1434.00 {(3033.56)} & 15.60 {(6.08)} & 1382.00 {(1200.85)} & 32.80 {(20.15)} & 3259.40 {(4395.80)} & 37.50 {(21.95)} \\
		Optimality gap & {0\%} & {0\%} & {0\%} & {0\%} & {0\%} & {0\%}
		\\
		\bottomrule
\end{tabular}}
\caption{Average results from solving models $\mathcal{M}$ and $\mathcal{M}_{sym}$ on small instances. Standard deviations are reported in brackets.} \label{tab_summary_results_small_instances}
\end{table}

\begin{table}[h!]
\centering
\resizebox{0.75\textwidth}{!}{
	\begin{tabular}{l|ll|ll}\toprule
		Number of bins & 50 ($\mathcal{M}$) & 50 ($\mathcal{M}_{sym}$) & 121 ($\mathcal{M}$) & 121 ($\mathcal{M}_{sym}$) \\ \hline
		Binary variables & 13000  & 13510 & 74415 & 75635 \\
		Continuous variables & 164350  & 167450 & 930369 & 937871 \\
		Equality constraints & 6755  & 170986 & 16340 & 946164 \\
		Inequality constraints & 471200  & 161200 & 2738230 & 922746 \\ \hline
		Profit (\euro) & 501.29 & 581.59 & $-$ & $-$ \\
		Weight of collected waste (kg) & 2306.69 & 2585.16 & $-$ & $-$ \\
		Travelled distance (km) & 190.72 & 193.96 & $-$ & $-$ \\
		Ratio weight/distance (kg/km) & 12.09 & 13.33 & $-$ & $-$\\
		CPU time (s) & 86400 & 86400 & OOM & OOM\\
		Optimality gap & 20.82\% & 2.77\% & $-$ & $-$\\
		\bottomrule
\end{tabular}}
\caption{Average results from solving models $\mathcal{M}$ and $\mathcal{M}_{sym}$ on the large instance. OOM stands for ``Out-Of-Memory''.} \label{tab_summary_results_big_instance}
\end{table}

\subsection{The impact of uncertainty and the quality of the deterministic solution} \label{sec_stochastic_measure}
The purpose of this section is twofold. Firstly, we discuss the importance of including stochasticity in the waste collection problem under investigation by comparing the optimal value of the stochastic formulation, i.e. the \emph{Recourse Problem} ($RP$) with the perfect information case (the so-called \emph{Wait and See} approach, $WS$, see \cite{BirLou2011}). Secondly, we show the benefits of considering stochasticity in model $\mathcal{M}$ with respect to its deterministic counterpart (the so-called \emph{Expected Value} problem, $EV$, see \cite{BirLou2011}. All the results reported below are obtained with small instances, as they are solved at optimality.

In the perfect information case, the realization of the waste accumulation rate is already known at the very first stage. Therefore, it is possible to compute the $WS$ as the average of the optimal values over single-scenario problems. The comparison between the $WS$ and the $RP$ is provided by the \emph{\%Expected Value of Perfect Information} ($\%EVPI$, see \cite{BirLou2011}), computed as
\begin{linenomath}
\begin{equation*}
	\%EVPI := (WS-RP)/RP.
\end{equation*}
\end{linenomath}

The results reported in the second column of Table \ref{tab_summary_small_instances} show that, on average, the $EVPI$ is 81\% of the $RP$. This means that, for reaching perfect information on the accumulation rate, the decision maker would be ready to pay at most 81\% of the total profit. Detailed results for each instance are shown in the Appendix.

\begin{table}[h!]
\centering
\resizebox{\textwidth}{!}{
	\begin{tabular}{c|c|c|c|ccc}\toprule
		Size & $\%EVPI$ & $\%VSS^t$, $1\leq t \leq 5$ & $\%MLUSS^t$, $1\leq t \leq 5$ & $\%MLUDS^1$ & $\%MLUDS^t$, $2\leq t \leq 4$ & $\%MLUDS^5$ \\ \hline
		9 bins & 188\% & $\infty$ & $\infty$ & 8\% & 674\% & 681\%\\
		10 bins & 36\% & $\infty$ & $\infty$ & 0\% & 175\% & 175\%\\
		11 bins & 20\% & $\infty$ & $\infty$ & 0\% & 158\% & 158\%\\ \hline
		Average & 81\% & $\infty$ & $\infty$ & 3\% & 336\% & 338\%\\
		\bottomrule
\end{tabular}}
\caption{Summary results of stochastic measures $\%EVPI$, $\%VSS^t$, $\%MLUSS^t$, $\%MLUDS^t$, for $1\leq t \leq 5$, expressed in percentage gap to the corresponding $RP$ problem.} \label{tab_summary_small_instances}
\end{table}

As a simpler approach, the decision maker may replace the waste accumulation rate by its expected value at each stage, and solve the deterministic $EV$ program. In a multi-stage context, the \emph{\%Value of Stochastic Solution at stage t} ($\%VSS^t$, see \cite{MagAllBer2014}), with $t\in\mathcal{T}'$, measures the expected gain from solving the stochastic model $RP$ rather than its deterministic counterpart up to stage $t$ as
\begin{linenomath}
\begin{equation*}
	\%VSS^t := (RP - EEV^t)/RP, \quad t = 1,\ldots,T-1.
\end{equation*}
\end{linenomath}

$EEV^t$ is the \emph{Expected result of using the EV solution until stage t} and denotes the objective function value of the $RP$ model, having fixed the decision variables $x_{ij}^t$ and $y_{i}^t$ on the routing until stage $t$ at the optimal values obtained by solving the $EV$ problem. In the great majority of the instances, the $EEV^t$ problems are infeasible and thus the corresponding $\%VSS^t$ is infinite already at the first stage (see column 3 of Table \ref{tab_summary_small_instances}). Indeed, by taking the average of the accumulation rate, the solution policy of the $EV$ problem may not impose any visiting schedules. On the other hand, the $RP$ model may require such a collection, being the accumulation rate high for certain bins: this is clearly a contradiction, leading to the infeasibility of the corresponding $EEV^t$ problem.

In the following, we further investigate if the deterministic solution still carries useful information for the stochastic case. To achieve this purpose, firstly we compute the \emph{Multi-stage Expected Skeleton Solution Value at stage t} ($MESSV^t$, see \cite{MagAllBer2014}), as the optimal objective function value of the $RP$ model having fixed at zero all the routing variables $x_{ij}^t$ and $y_i^t$ that are zero in the $EV$ problem until stage $t$. This allows to test whether the deterministic model provides the correct non-zero variables. Once the $MESSV^t$ is computed, it is compared with the $RP$ by introducing the \emph{$\%$Multi-stage Loss Using Skeleton Solution until stage t} ($\%MLUSS^t$), expressed as
\begin{linenomath}
\begin{equation*}
	\%MLUSS^t := (RP - MESSV^t)/RP, \quad t = 1,\ldots,T-1.
\end{equation*}
\end{linenomath}

The results in Table \ref{tab_summary_small_instances} and in the Appendix show that $\%VSS^t$ coincides with $\%MLUSS^t$ for all $t=1,\ldots,5$, both in the case of infiniteness and finiteness of $\%VSS^t$. On the one hand, $\%MLUSS^t=\%VSS^t=\infty$ if the average data on the waste accumulation rate do not support a collection in the $EV$ solution: $EEV^t$ and $MESSV^t$ problems are both infeasible. On the other hand, if $\%VSS^t$ is finite, the deterministic model correctly selects the bins to be visited and their combination into routing plans. Therefore, for all the non-visited bins, the values of the decision variables $x_{ij}^t$ and $y_i^t$ are fixed to zero in the $MESSV^t$ problem. This implies straightforwardly that there is only one possible path for the vehicle to visit the selected bins.

Finally, we carry out an analysis regarding the upgradeability of the expected value solution to become good, or optimal, in the stochastic setting. Specifically, we consider the $EV$ solution $\bar{x}_{ij}^t$, $\bar{y}_{i}^t$ until stage $t$ as a starting point in the $RP$ model, by adding the constraints $x_{ij}^t \geq \bar{x}_{ij}^t$, for all $i,j\in\mathcal{I}$, and $y_i^t \geq \bar{y}_{i}^t$, for all $i\in\mathcal{I}'$ up to stage $t$. The corresponding optimal value is denoted as \emph{Multi-stage Expected Input Value until stage t} ($MEIV^t$, see \cite{MagAllBer2014}). From this measure, the \emph{$\%$Multi-stage Loss of Upgrading the Deterministic Solution until stage t} ($\%MLUDS^t$) is defined as follows
\begin{linenomath}
\begin{equation*}
	\%MLUDS^t := (RP - MEIV^t)/RP, \quad t = 1,\ldots,T-1.
\end{equation*}
\end{linenomath}

As it is reported in Table \ref{tab_summary_small_instances}, $\%MLUDS^1$ is close to zero on average. Indeed, only in \emph{inst\_4\_9} (see the Appendix) $\%MLUDS^1$ is strictly positive. This situation derives from a collection after stage 2 in the $EV$ solution, where the conditions $x_{ij}^1 \geq 0$, for all $i,j\in\mathcal{I}$, and $y_{i}^1 \geq 0$, for all $i\in\mathcal{I}'$, are automatically satisfied by constraints \eqref{constraint_14}-\eqref{constraint_15} in the $MEIV^1$ problem. In all the other instances, at stage 2 the $EV$ problem imposes a collection on a subset of bins with respect to the $RP$ problem and, thus constraints \eqref{constraint_15} are themselves satisfied in the $MEIV^1$ problem.

The large values of $\%MLUDS^t$, with $t=2,\ldots,5$, depend on the fact that the corresponding $MEIV^t$ problems perform collections at the consecutive stages 1 and 2, due to the additional constraints on the $EV$ solution at stage 2. Being the amount of waste in the bins low since already emptied at stage 1, the transportation costs are greater than the revenues: $MEIV^t$ turns to be negative and hence $\%MLUDS^t$ reaches high values (on average, 338\%).

The results discussed so far justify the adoption of a stochastic model over a deterministic formulation when addressing a waste collection problem. From the analysis, we conclude that it is possible to take the deterministic solution as input in the stochastic model only in the first stage, whereas in the next stages the $EV$ solution is no longer to be upgradeable.

\subsection{Performance of the rolling horizon approach} \label{sec_performance_RH}
Since model $\mathcal{M}$ is NP-hard and with large instances obtaining the optimal solution is challenging (see Table \ref{tab_summary_results_big_instance}), in this section we evaluate the performance of the rolling horizon approach both in terms of average profit reduction and CPU time savings. Instead of solving a $T$-stage stochastic program, such heuristic considers a sequence of $T-1$ subproblems defined over a reduced number of stages. In our case study, $T=6$ and, thus the reduced number of periods $W$ is an integer between 1 and 4. Detailed results for each small instance are presented in the Appendix.

\begin{figure}[h!]
\centering
\includegraphics[width=0.8\textwidth]{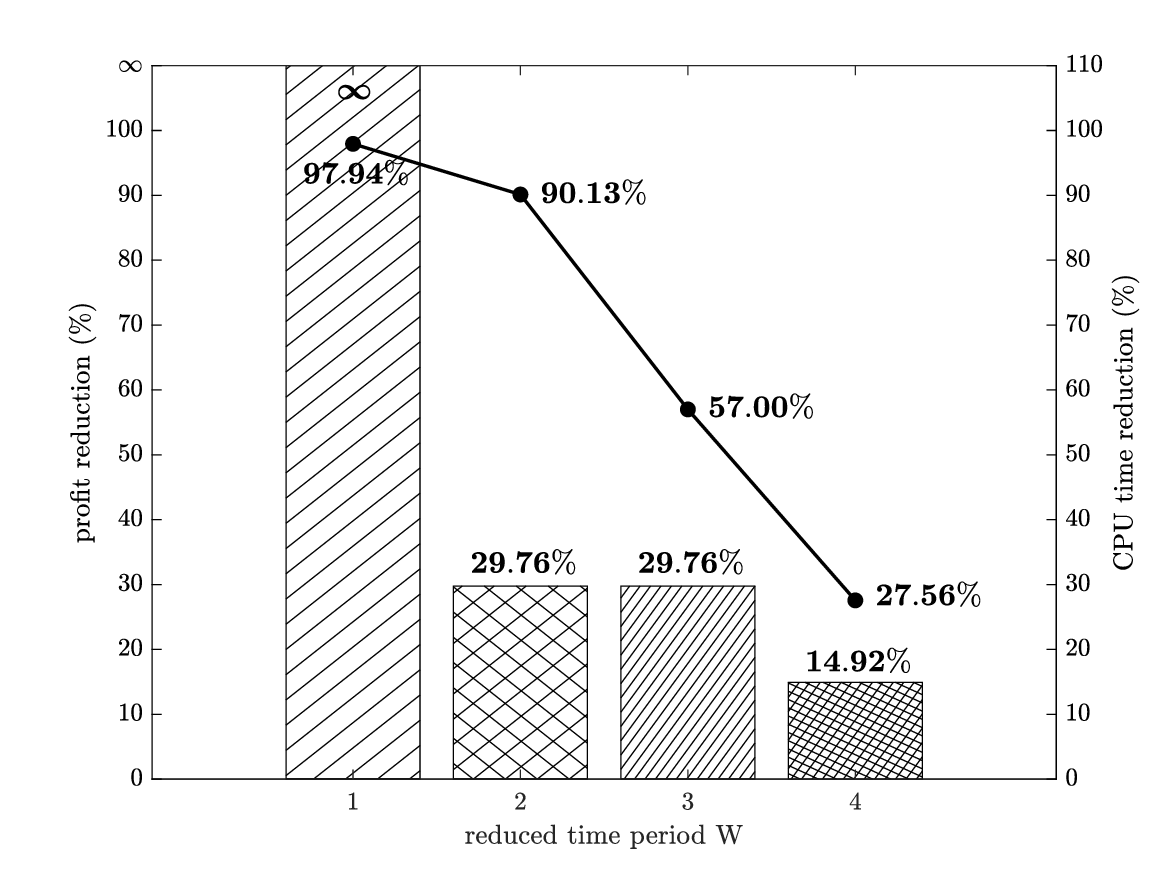}
\caption{Performance of the rolling horizon approach for the small-sized instances. The vertical bars represent the profit percentage reduction when applying the rolling horizon approach (left-hand scale). The results show the average over the thirty instances. When $W=1$, due to infeasibility, the reduction may be infinite. The solid line refers to the CPU time percentage reduction to solve at optimality with the rolling horizon approach, compared to the original six-stage program (right-hand scale).}
\label{fig_RH_red_profit_time}
\end{figure}

In Figure \ref{fig_RH_red_profit_time} (left-hand scale) we depict with vertical bars the average percentage gap between the optimal $RP$ objective function value and the heuristic solution. As highlighted in Section \ref{sec_rolling_horizon}, when $W=1$ the rolling horizon approach may be infinitely suboptimal. Indeed, over the thirty small instances, five of them exhibit infeasibility in the first two-stage problem. Furthermore, the results obtained with $W=2$ and $W=3$ (29.76\% in both cases) are close to the one obtained for $W=1$ in the twenty-five feasible instances (29.53\%). However, with $W=2$ and $W=3$ no infeasibility issues occur in any instance. Finally, the performance of the heuristic improves when $W=4$, as the average profit gap is of 14.92\%.

Regarding the computational time, we report in Figure \ref{fig_RH_red_profit_time} (right-hand scale) the results in terms of average percentage reduction with respect to the six-stage model. We notice that the CPU time reduction decreases, as $W$ goes from 1 to 4, due to the increasing size of the subproblems. Specifically, when $W=1$ and $W=2$ the average savings are of 97.94\% and 90.13\%, respectively. The CPU time reduction is of 57.00\% with $W=3$, even if the performance in terms of profit is the same as $W=2$. Lastly, the similarity of the optimal profit between the $RP$ model and the rolling horizon approach with $W=4$ requires a significant computational effort, as the average CPU time saving is of 27.56\%.

From the previous analysis, we conclude that the rolling horizon approach is effective for the proposed multi-stage stochastic model. As expected, the performance of such heuristic strongly depends on the size $W$ of the reduced time horizon. If the decision maker requires a good accuracy in a short time, $W=2$ is the best candidate. On the other hand, if she/he is willing to wait, $W=4$ attains better results but in a longer computational time.

\subsection{A real case study} \label{sec_real_case_study}
In this section, we present the results of the simulations in a real case study. We consider a large instance composed by 50 bins randomly chosen from the original set of 121 waste containers.

Given the high number of variables and of constraints (see Table \ref{tab_summary_results_big_instance}), Gurobi is not able to solve at optimality model $\mathcal{M}$ within the time limit of 86400 seconds (one day), providing a relative optimality gap of 20.82\%. However, when considering model $\mathcal{M}_{sym}$ with distance the average between $d_{ij}$ and $d_{ji}$ the results are considerably better: the optimality gap is reduced to 2.77\% (see Table \ref{tab_summary_results_big_instance}).

For this reason, we deal with model $\mathcal{M}_{sym}$ for the study of the large case instance. Specifically, we investigate the performance of the rolling horizon approach under a reduced runtime limit, as one day may be excessively high from a managerial perspective. Thus, we set a time limit $TL$ of 2, 4, 6, 12, 24 hours on the whole algorithm and, accordingly, a time limit $TL_{sub}$ for each of the corresponding subproblems, as $TL_{sub}=\frac{TL}{\# \text{ }subproblems}$. Following the approach of \cite{CavBerMag2022}, after solving a subproblem, if some time is left, we add the remaining time to the following subproblem to be solved.

\begin{figure}[h!]
\centering
\includegraphics[width=0.8\textwidth]{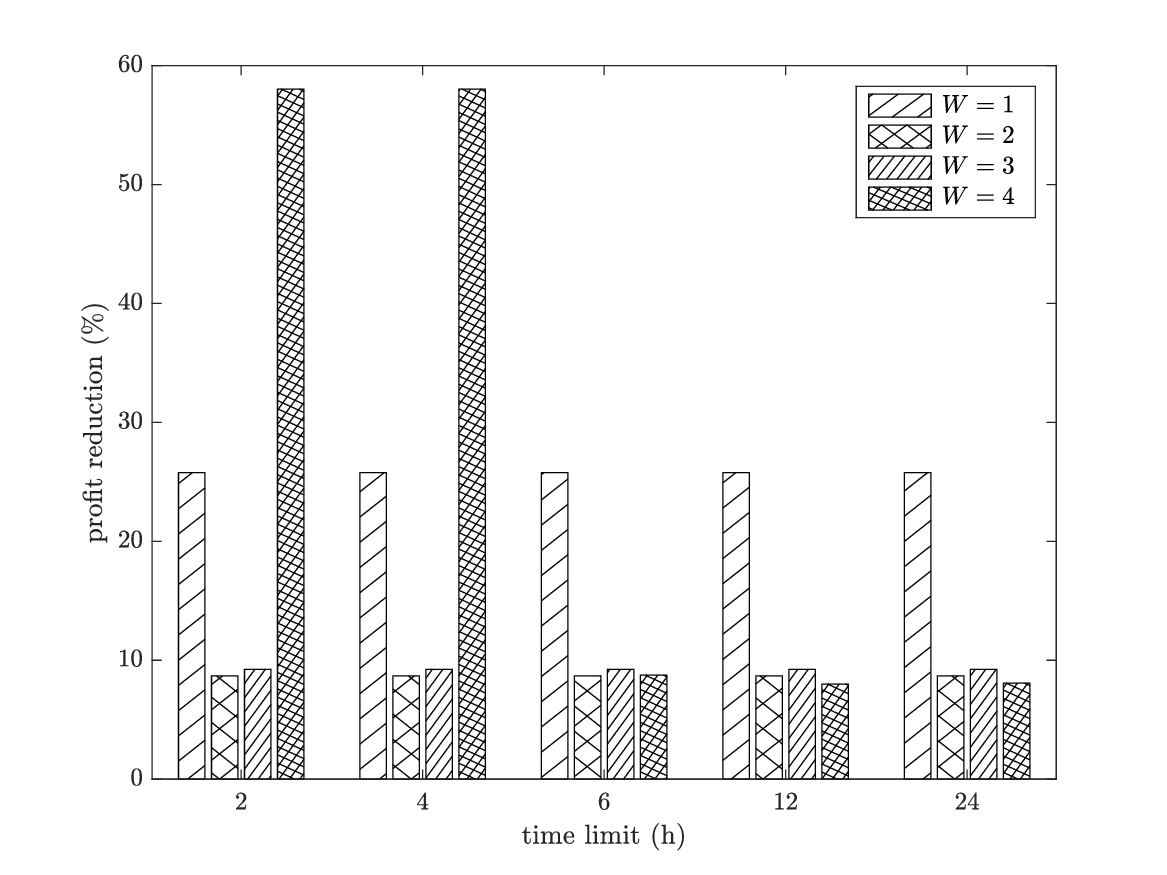}
\caption{Performance of the rolling horizon approach for the instance with 50 bins in terms of profit reduction.}
\label{fig_RH50_red_profit}
\end{figure}

Figure \ref{fig_RH50_red_profit} displays the percentage profit reduction, when applying the rolling horizon approach over a reduced time period $W$ and with different runtime limits $TL$. The reduction is with respect to the objective function value obtained when solving the six-stage model $\mathcal{M}_{sym}$ in one day. On the one hand, the rolling horizon approach with $W=1,2,3$ does not improve its performance when increasing the time limit. Specifically, when $W=1$ the profit reduction is 25.77\%, while for $W=2$ and $W=3$ it is of 8.68\% and 9.24\%, respectively, regardless of the runtime limit. On the other hand, when the time limit is enlarged, the heuristic with $W=4$ shows an enhancement of the results: from a reduction of 58.02\%, when $TL$ is 2 and 4 hours, to 8.07\% with $TL$ equal to 12 and 24 hours. The bad performance with low runtime limits (2 and 4 hours) is due to the large size of the first two subproblems, defined respectively on stages 1-5 and 2-6, which are challenging to solve in a short time ($TL_{sub}$ is equal to 24 and 48 minutes, respectively).

\begin{figure}[h]
\centering
\includegraphics[width=0.8\textwidth]{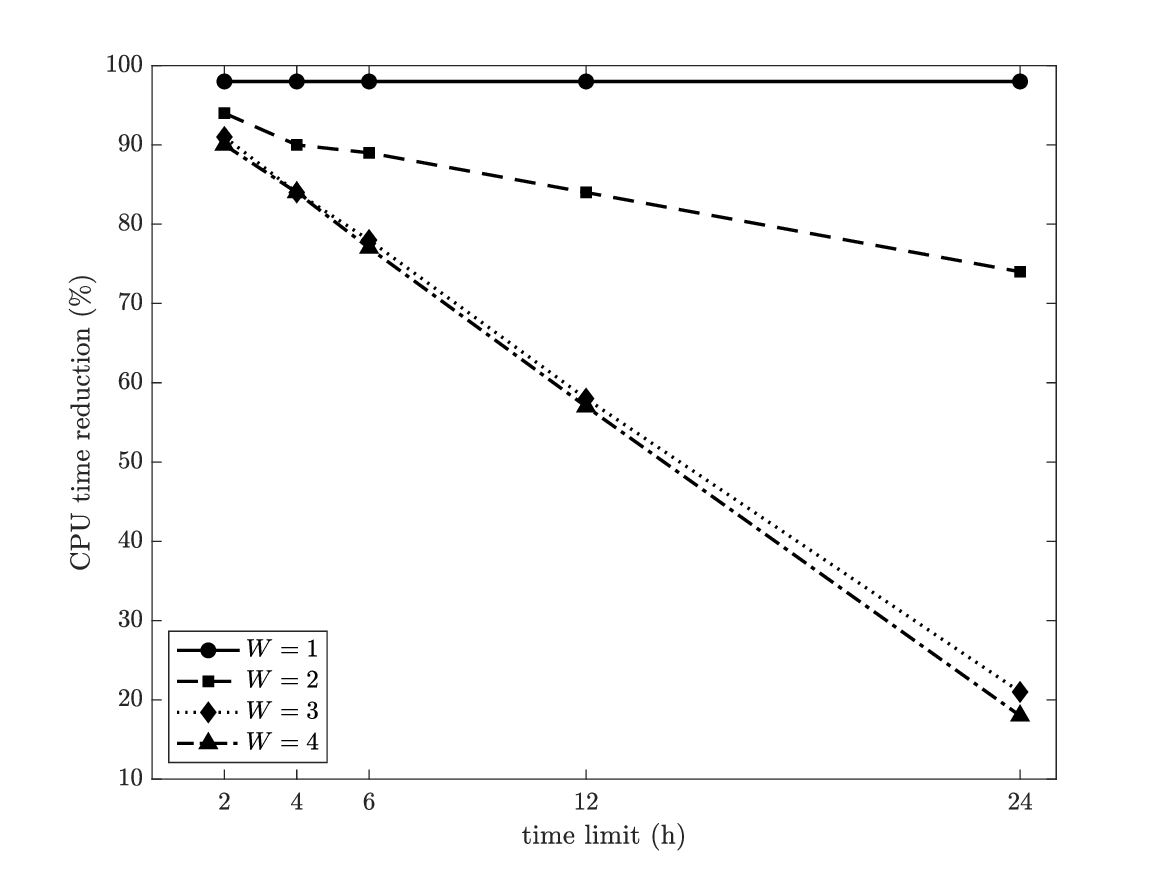}
\caption{Performance of the rolling horizon approach for the instance with 50 bins in terms of CPU time reduction.}
\label{fig_RH50_red_time}
\end{figure}

Similarly to the analysis carried out on the small instances, in Figure \ref{fig_RH50_red_time} we depict the percentage CPU time saving . When $W=1$, the reduction is equal to 98.34\%, independently of the runtime limit. Indeed, all the five subproblems are solved within the time limit $TL_{sub}$. If $W=2$, the CPU time saving is high with $TL$ equal to 2 hours (93.78\%), while it reaches a minimum of 73.91\% with $TL=24$ hours. Finally, the situations with $W=3$ and $W=4$ show a similar linear trend. Indeed, in these cases, the runtime limit is always reached, because of the large size of the subproblems.

\begin{figure}[h!]
\centering
\begin{subfigure}[b]{0.49\textwidth}
	\centering
	\includegraphics[width=\textwidth]{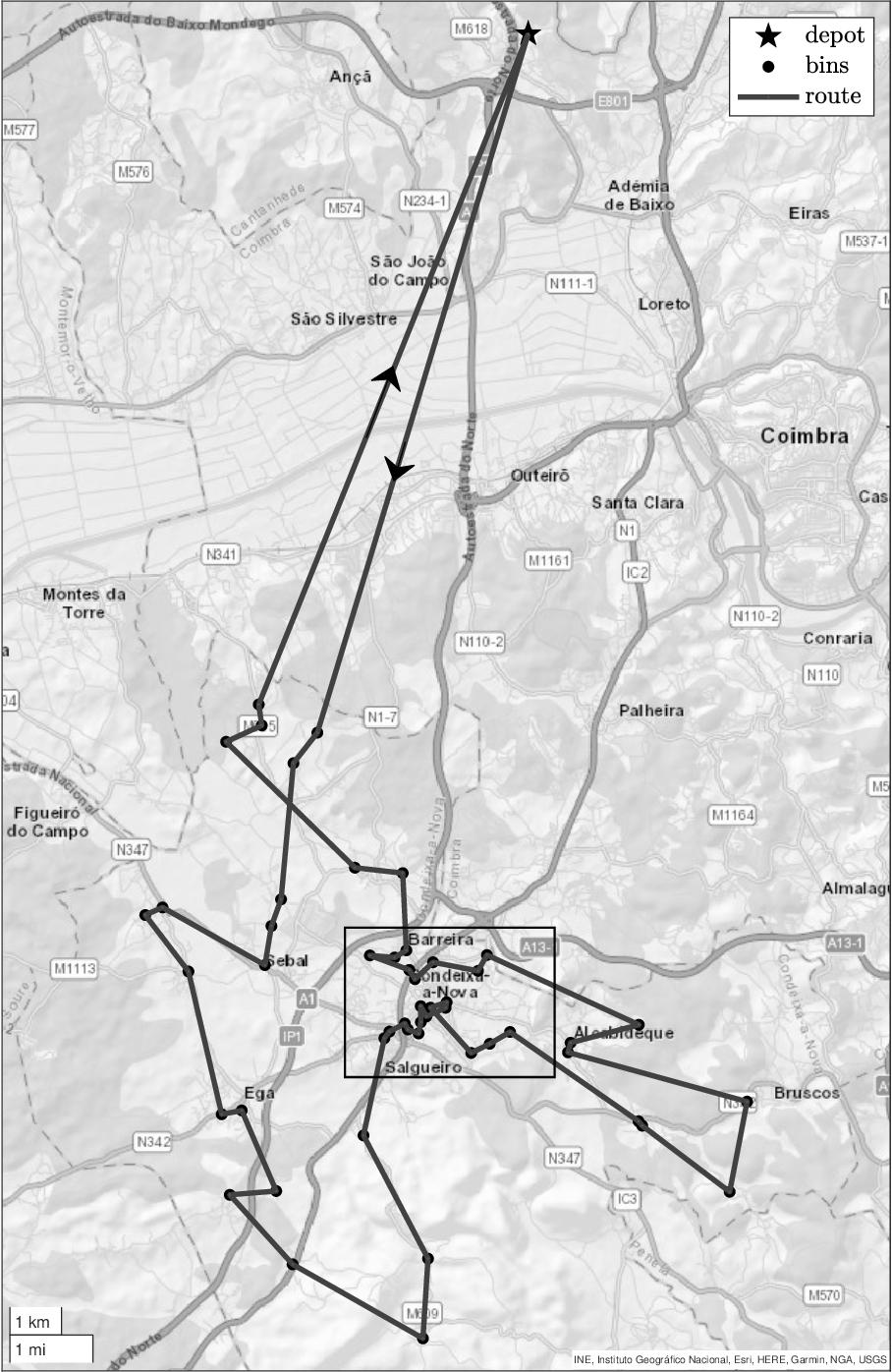}
\end{subfigure} \hfill
\begin{subfigure}[b]{0.49\textwidth}
	\centering
	\includegraphics[width=\textwidth]{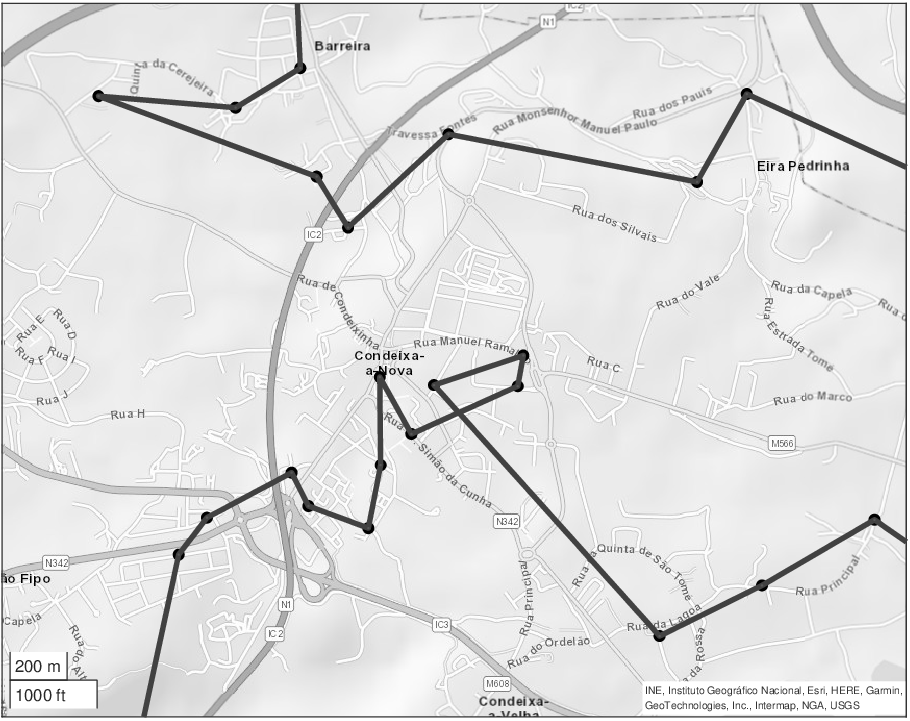}
\end{subfigure}

\caption{Route for the large case instance with 50 bins, obtained by applying the rolling horizon approach with $W=2$ and runtime limit of 2 hours. The route is performed on days 2 and 6 of the planning period. In the picture on the right, a zoom on the area of \textit{Condeixa-a-Nova} is depicted.}
\label{fig_RH50_route_W2}
\end{figure}

By combining all the previous observations, we conclude that it is worth applying the rolling horizon approach when solving large instances of the proposed multi-stage stochastic waste collection problem. The performance depends not only on the reduced time horizon of such heuristic but also on the runtime limit set by the user. Confirming the results obtained with small instances, we conclude that the rolling horizon approach with $W=2$ and runtime limit of 2 hours is a good trade-off between accuracy and time savings. For the sake of illustration, in Figure \ref{fig_RH50_route_W2} we depict the route obtained with these choices of the parameters for the large case instance.

\subsection{Managerial insights} \label{sec_managerial_insights}
We conclude our analysis by providing some managerial insights on the discussed problem.

In Table \ref{tab_managerial_insights} we report key performance indicators obtained from the proposed formulation when considering a runtime limit of two hours for the rolling horizon approach. We notice that the highest value of both the profit (see the second row) and the ratio between the total weight of collected waste and the total travelled distance (see the fifth row) is reached in the case of $W=2$. This implies that choosing such reduced time period leads to a more efficient and cost-effective planning when compared to the other cases. A similar conclusion can be drawn when considering the ratio between the profit reduction and the CPU time reduction with respect to the original $RP$ problem (see the sixth row).

\begin{table}[h!]
\centering
\resizebox{0.9\textwidth}{!}{
	\begin{tabular}{lllll}\toprule
		$W$ & 1 & 2 & 3 & 4\\ \hline
		Profit (\euro) & 431.70	& {531.13} & 527.82 & 244.17\\
		Total weight of collected waste (kg) & 2573.31 & 2558.26 & 2564.08 & 1201.77\\
		Total travelled distance (km) & 340.29 & 236.35 & 241.41 & 116.36\\
		Weight/distance (kg/km) & 7.56 & {10.82} & 10.62 & 10.33\\
		Profit reduction/CPU time reduction & 0.26 & 0.09 & 0.10 & 0.64\\
		\bottomrule
\end{tabular}}
\caption{Key performance indicators for the real case instance of 50 bins, when applying the rolling horizon approach with a runtime limit of 2 hours.} \label{tab_managerial_insights}
\end{table}

\begin{figure}[h!]
\centering
\begin{subfigure}[b]{0.49\textwidth}
	\centering
	\includegraphics[width=\textwidth]{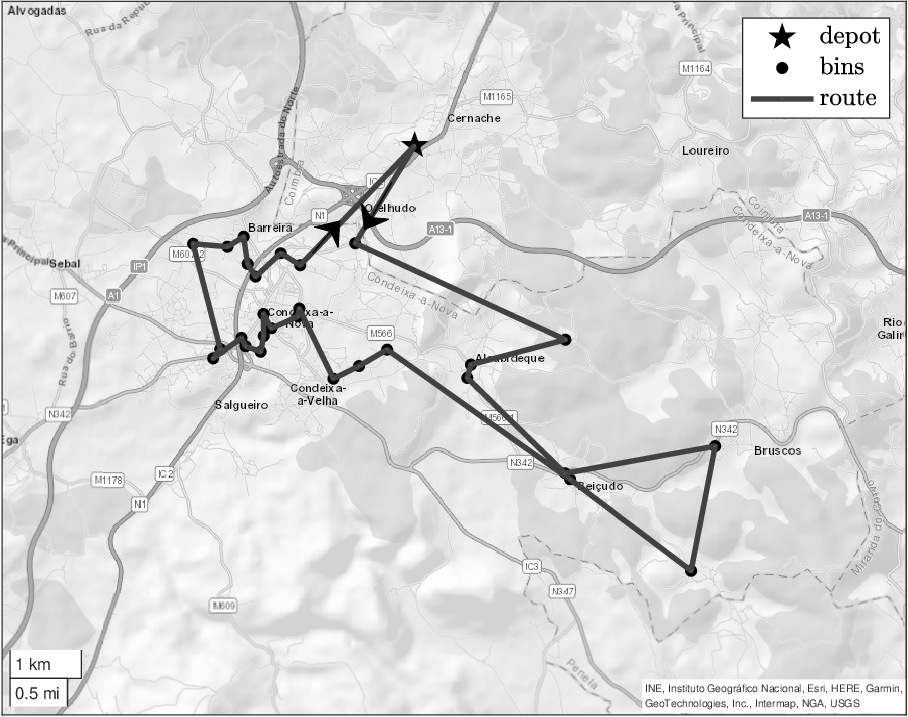}
	\caption{Day 2}
\end{subfigure} \hfill
\begin{subfigure}[b]{0.49\textwidth}
	\centering
	\includegraphics[width=\textwidth]{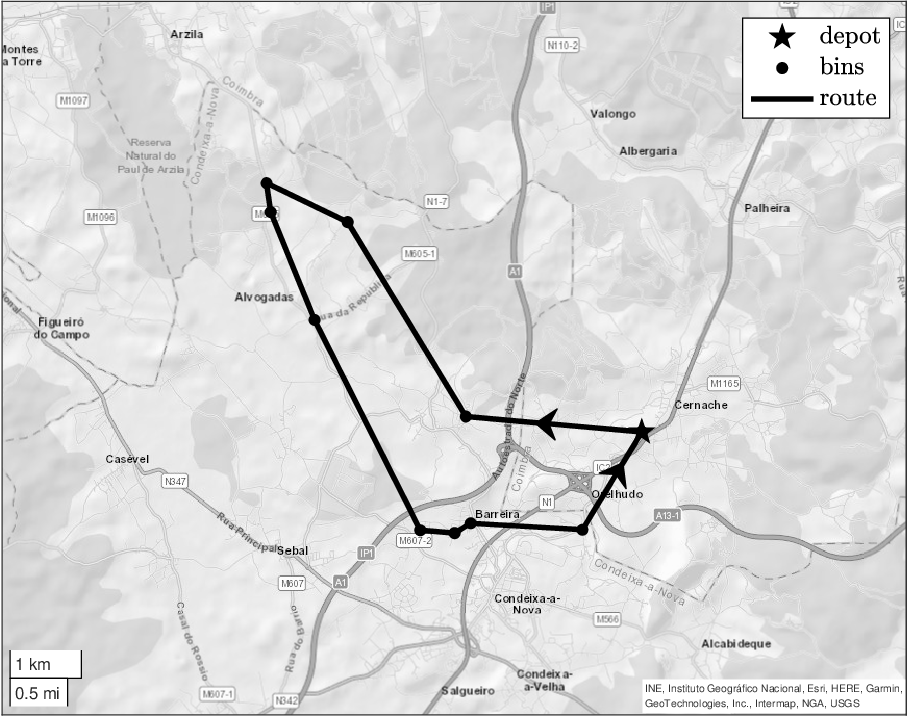}
	\caption{Day 5}
\end{subfigure}
\begin{subfigure}[b]{0.49\textwidth}
	\centering
	\includegraphics[width=\textwidth]{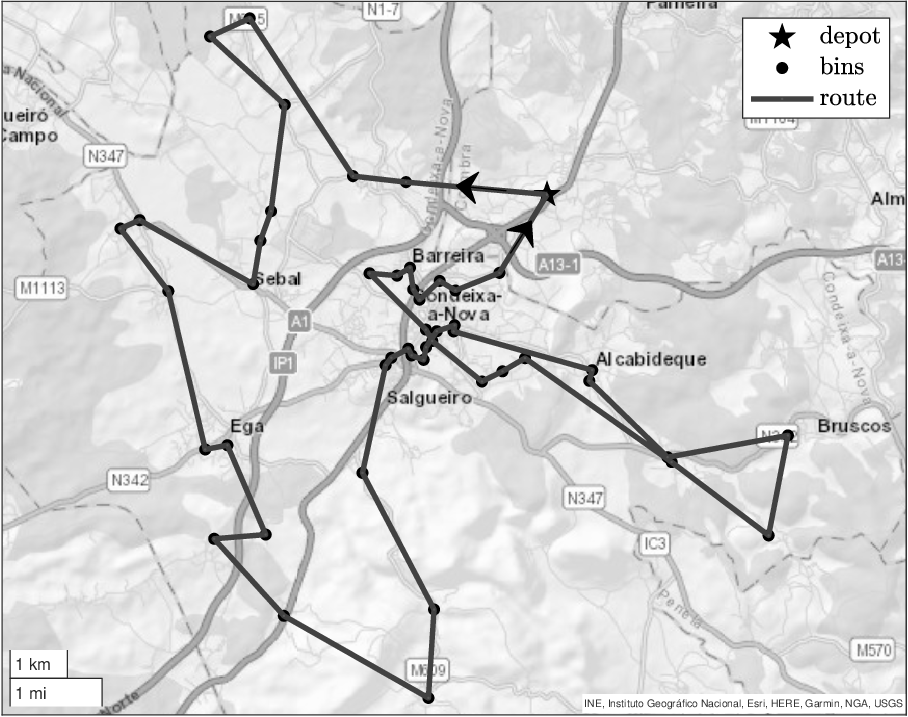}
	\caption{Day 6}
\end{subfigure} \hfill
\caption{Routes performed on days 2, 5, 6 of the planning period, with a closer depot to the fifty bins. The results are obtained by applying the rolling horizon approach with $W=2$ and runtime limit of 2 hours.}
\label{fig_RH50_route_W2_centreddepot}
\end{figure}

From a managerial perspective, one of the key feature of the proposed model is the selection of the bins to be visited. In Figure \ref{fig_RH50_route_W2} the same route is performed on days 2 and 6 of the planning period. Almost all the bins are visited twice (50 and 49 bins, respectively, on day 2 and 6), due to the very high distance between the depot and the considered municipality. In Figure \ref{fig_RH50_route_W2_centreddepot} we depict the results of a simulation with the same bins as in the instance described in Section \ref{sec_real_case_study}, but with a closer depot. We notice that the collection is now performed on three days, with an accurate selection of bins, respectively 28, 9 and 47 bins for each collection day. The profit is increased by 15.26\% (see the fourth column of Table \ref{tab_managerial_insights_closerdepot_varuncert}), due to the strong decrease of the total travelled distance (151.62 km vs. 236.35 km). The total weight of collected waste remains almost unchanged (2545.94 kg vs. 2558.26 kg), implying an increase of the ratio \emph{weight over distance} by 55.18\%. The waste manager benefits from these results since they suggest that the opening of a new depot, closer to the bins, increases significantly the profit, with more selective and accurate routes.

Inspired by the work of \cite{ElbWoh2016}, we performed a sensitivity analysis with respect to the bins' filling level. Specifically, we considered situations where the accumulation rate is increased or decreased by 10\%. This corresponds to a variation of the uncertainty in the construction of the scenario trees. The results of the simulations are reported in the last columns of Table \ref{tab_managerial_insights_closerdepot_varuncert}. With respect to the original setting, the increase of the amount of waste in the bins implies a 7.10\% higher profit for the waste company. On the other hand, when the amount of waste is reduced, the total profit decreases by 8.55\%. Similar considerations can be formulated for all the other key performance indicators.

\begin{table}[h!]
\centering
\resizebox{\textwidth}{!}{
	\begin{tabular}{ll|lc|lc|lc}\toprule
		& {Original setting} & \multicolumn{2}{c|}{{Closer depot}} & \multicolumn{4}{c}{{Varying bins' filling level}}\\
		& & & {$\Delta \%$} & {$+10\%$} & {$\Delta \%$} & {$-10\%$} & {$\Delta \%$}
		\\ \hline
		{Profit (\euro)} & {531.13} & {612.17} & {$+15.26\%$} & {568.85} & {$+7.10\%$} & {485.70} & {$-8.55\%$}
		\\
		Total weight of collected waste (kg) & {2558.26} & {2545.94} & {$-4.82\%$} & {2696.10} & {$+5.39\%$} & {2406.82} & {$-5.92\%$}
		\\
		{Total travelled distance (km)} &  {236.35}  & {151.62} & {$-35.85\%$} & {239.98} & {$+1.54\%$} & 236.35 & {$0\%$}
		\\
		{Weight/distance (kg/km)} & {10.82} & {16.79} & {$+55.18\%$} & {11.23} & {$+3.79\%$} & {10.18} & {$-5.91\%$}
		\\
		\bottomrule
\end{tabular}}
\caption{Key performance indicators with different configurations of the 50 bins instance. All the results are obtained applying the rolling horizon approach with $W=2$ and a runtime limit of 2 hours. The percentage variation $\Delta \%$ is computed with respect to the the original setting in the second column.} \label{tab_managerial_insights_closerdepot_varuncert}
\end{table}

\section{Conclusions} \label{sec_conclusions}
In this paper, we have studied a Stochastic Inventory Routing Problem for waste collection of recyclable materials with uncertain bins' filling level. We have formulated such a kind of problem through a multi-stage mixed-integer stochastic programming linear model, with the aim of maximizing the total expected profit. Scenario trees on waste accumulation rate have been generated by means of conditional density estimation and dynamic stochastic approximation techniques. A validation in terms of in-sample stability has been assessed too. The impact of stochasticity in the proposed waste collection problem has been investigated through standard stochastic measures, showing the benefits of the stochastic formulation when compared to the deterministic case. To face with the computational complexity of the problem, we have proposed the rolling horizon as heuristic methodology and derived a worst-case analysis of its performance. Given the availability of real data, we have carried out extensive computational experiments on small- and large-sized instances. We have tested the performance of the rolling horizon approach, founding out that, if the reduced time horizon is properly chosen, such heuristic provides good quality solutions with limited computational efforts. Finally, we have drawn managerial insights considering different configurations in a real case instance and providing key performance indicators.

Regarding future developments, several streams of research can originate from this work. First of all, given the recent availability of data on a near-continuous basis, a stochastic programming model with real-time information provided by sensors should be formulated. Secondly, the size of the test instances should be enlarged to include more realistic waste collection systems. As a downside, this would make the proposed formulation very challenging from a computational perspective. To this extent, Benders' decomposition and column generation algorithms would be useful techniques to be considered in the future.

\section*{Data}
The data used in this article are publicly available on a GitHub repository (\url{https://github.com/aspinellibg/StochwasteIRP}).

\section*{Acknowledgments}
Research activities of Andrea Spinelli were performed both at University of Bergamo (Italy) and at Instituto Superior T\'ecnico (Portugal), thanks to the international mobility grant of University of Pavia, year 2022.

This work has been supported by ``ULTRA OPTYMAL - Urban Logistics	and sustainable TRAnsportation: OPtimization under uncertainTY	and MAchine Learning'', a PRIN2020 project funded by the Italian	University and Research Ministry (grant no. 20207C8T9M, official website: https://ultraoptymal.unibg.it), by the Portuguese Foundation for Science and Technology (FCT) through the research projects no. 2022.04180.PTDC and UIDB/00097/2020, by European Union Next Generation EU - PRIN2022 (grant no. 20223MHHA8), by H2020 TUPLES, and by USAF AFORS (grant no. FA8655-21-1-7046).

\section*{References}
\bibliography{StochasticWaste_bib}

\newpage
\section*{Appendix}

\emph{A)} \emph{Multi-stage stochastic model $\mathcal{M}_{sym}$ with a two-commodity flow formulation}\\
\emph{\underline{Sets}:}\\
$\mathcal{I}=\{i:i=0,1,\ldots,N,N+1\}$: set of $N$ waste bins, the real depot 0 and the copy depot $N+1$;\\
$\mathcal{I'}=\{i:i=1,\ldots,N\}$: set of $N$ waste bins (depots excluded);\\
$\mathcal{T}=\{t:t=1,\ldots,T\}$: set of stages;\\
$\mathcal{T'}=\{t:t=1,\ldots,T-1\}$: set of stages (last stage excluded);\\
$\mathcal{T''}=\{t:t=2,\ldots,T\}$: set of stages (first stage excluded);\\
$\mathcal{N}^1=\{n:n=1\}$: root node at stage 1;\\
$\mathcal{N}^t=\{n:n=1,\ldots,n^t\}$: set of ordered nodes of the tree at stage $t\in\mathcal{T}$.\\
\emph{\underline{Deterministic parameters}:}\\
$C$: travelling cost per distance unit;\\
$R$: selling price of a recyclable material;\\
$Q$: vehicle capacity;\\
$B$: waste density;\\
$M$: Big-M number;\\
$d_{ij}$: distance between $i\in\mathcal{I}$ and $j\in\mathcal{I}$;\\
$S_{i}^{init}$: percentage of waste on the total volume of bin $i\in\mathcal{I'}$ at the first stage;\\
$E_i$: capacity of bin $i\in\mathcal{I'}$;\\
$pa(n)$: parent of node $n\in\mathcal{N}^t$, $t\in\mathcal{T''}$.\\
\emph{\underline{Stochastic parameters}:}\\
$a_i^{n}$: uncertain accumulation rate of bin $i\in\mathcal{I'}$ at node $n\in\mathcal{N}^t$, $t\in\mathcal{T''}$;\\
$\pi^{n}$: probability of node $n\in\mathcal{N}^t$, $t\in\mathcal{T}$.\\
\emph{\underline{Decision variables}:}\\
$x_{ij}^t\in\{0,1\}$: binary variable indicating if arc $(i,j)$ is visited at time $t+1$, with $t\in\mathcal{T'}$ and for $i,j\in\mathcal{I}$, $i\neq j$;\\
$y_{i}^t\in\{0,1\}$: binary variable indicating if waste bin $i\in\mathcal{I'}$ is visited at time $t+1$, with $t \in\mathcal{T'}$;\\
$f_{ij}^{n}\in\mathbb{R}^{+}$: non-negative variable representing the waste flow between $i\in\mathcal{I'}$ and $j\in\mathcal{I}$, $i\neq j$, for $n\in\mathcal{N}^t$, $t\in\mathcal{T''}$;\\
$w_i^{n} \in \mathbb{R}^{+}$: non-negative variable representing the amount of waste collected at waste bin $i\in\mathcal{I'}$, for $n\in\mathcal{N}^t$, $t\in\mathcal{T''}$;\\
$u_i^{n} \in \mathbb{R}^{+}$: non-negative variable representing the amount of waste at waste bin $i\in\mathcal{I'}$, for $n\in\mathcal{N}^t$, $t\in\mathcal{T}$.

\emph{\underline{Model $\mathcal{M}_{sym}$}:}
\begin{eqnarray*} &	\max & R\sum_{t\in\mathcal{T''}} \sum_{n\in\mathcal{N}^t} \pi^{n} \sum_{i\in\mathcal{I'}} w_i^{n} -{\color{black}{\frac{C}{2}}}\sum_{t\in\mathcal{T'}}\sum_{\substack{i,j\in\mathcal{I}\\i\neq j}}d_{ij}x_{ij}^t \\
&	& \text{s.t. } \sum_{\substack{j\in\mathcal{I}\\j\neq i}} (f_{ij}^{n} - f_{ji}^{n}) = {\color{black}{2}}w_{i}^{n}\qquad i\in\mathcal{I'}, n\in\mathcal{N}^t, t\in\mathcal{T''}\\
&	& {\color{black}{\sum_{i\in\mathcal{I'}} f_{iN+1}^{n} = \sum_{i\in\mathcal{I'}}w_{i}^{n} \qquad n\in\mathcal{N}^t, t\in\mathcal{T''}}}\\
&	& {\color{black}{f_{ij}^{n}+f_{ji}^n = Qx_{ij}^{t-1} \qquad i,j \in \mathcal{I}, i\neq j, n\in\mathcal{N}^t, t\in\mathcal{T''}}}\\
&	& f_{ij}^{n} \leq (Q-E_jBa_{j}^{n})x_{ij}^{t-1} \qquad i,j \in \mathcal{I}', i\neq j, n\in\mathcal{N}^t, t\in\mathcal{T''}\\
&	& {\color{black}{\sum_{\substack{i\in\mathcal{I}\\i\neq j}} x_{ij}^t = 2y_j^t \qquad j\in\mathcal{I'}, t\in\mathcal{T'}}}\\
&	& w_i^{n} \leq E_iBy_{i}^{t-1} \qquad i\in\mathcal{I'}, n\in\mathcal{N}^t, t \in \mathcal{T''}\\
&	& u_i^{n} \leq M(1- y_{i}^{t-1}) \qquad i\in\mathcal{I'}, n\in\mathcal{N}^t, t \in \mathcal{T''}\\
&	& u_i^{n} = E_iBS_i^{init} \qquad i\in\mathcal{I'}, n\in\mathcal{N}^1\\
&	& u_i^{n} = u_i^{pa(n)}+E_iBa_i^{n}-w_i^{n} \qquad i\in\mathcal{I'}, n\in\mathcal{N}^t, t\in\mathcal{T''} \\
&	& u_i^{pa(n)} \leq \big(1 -a_i^{n}\big)E_iB \qquad i\in\mathcal{I'}, n\in\mathcal{N}^t, t\in\mathcal{T''}\\
&	& x_{ij}^t \in\{0,1\} \qquad i,j\in\mathcal{I}, i\neq j, t\in\mathcal{T'}\\
&	& y_{i}^t \in \{0,1\} \qquad i\in\mathcal{I'}, t\in\mathcal{T'}\\
&	& f_{ij}^{n} \geq 0 \qquad i\in\mathcal{I}', j\in\mathcal{I}, i\neq j, n\in\mathcal{N}^t, t\in\mathcal{T''}\\
&	& w_{i}^{n} \geq 0 \qquad i\in\mathcal{I'}, n\in\mathcal{N}^t, t\in\mathcal{T''}\\
&	& u_i^{n} \geq 0 \qquad i\in\mathcal{I'}, n\in\mathcal{N}^t, t\in\mathcal{T}
\end{eqnarray*}
\clearpage
\emph{B)} \emph{Scenario tree generation}

In this section, we discuss how to generate scenario trees to describe the problem uncertainty. We adopt the methodologies proposed in \cite{KirPflPic2020}, which are based on the works of Pflug and Pichler (see \cite{PflPic2016} for futher details).

Since only a limited number of trajectories of the accumulation rate is available from historical data, new and additional samples are needed to be generated, even if the true distribution of the accumulation rate is not known. However, it can be estimated by a non-parametric kernel density technique discussed in the following.

Let $\big(a_{i,o}^{(1)},\ldots,a_{i,o}^{(t)},\ldots,a_{i,o}^{(T)}\big)$ be the vector denoting the accumulation rate of bin $i$ for week of observation $o$, with $o=1,\ldots,N_o$. Let $k(\cdot)$ be a kernel function and $\big(p_1,\ldots,p_o,\ldots,p_{N_o}\big)$ be a $N_o$-dimensional vector of positive weights such that $\sum_{o=1}^{N_{o}}p_o=1$. Let $\alpha$ be a random number drawn from the uniform distribution $\mathcal{U}\big(0,1\big)$. At stage $t=1,\ldots,T$, a new sample $\widehat{a}_{i}^{(t)}$ of the accumulation rate of bin $i$ is given by
\begin{linenomath}
\begin{equation*}
	\widehat{a}_{i}^{(t)} = a_{i,o^*}^{(t)}+h^{(t)} \cdot K^{(t)},
\end{equation*}
\end{linenomath}
where:
\begin{itemize}
\item $o^*$ is an index between $1$ and $N_o$ such that $ \sum_{o=1}^{o^*-1} p_o < \alpha \leq  \sum_{o=1}^{o^*}p_o$;
\item $h^{(t)}$ is the bandwidth, computed according to the Silverman's rule of thumb (see \cite{Sil98}), namely $h^{(t)}=\sigma^{(t)} \cdot N_{o}^{-\frac{1}{m^{(t)}+4}}$, being $\sigma^{(t)}$ the standard deviation of data at stage $t$ and $m^{(t)}$ the dimension of the process at stage $t$;
\item $K^{(t)}$ is a random value sampled from the kernel distribution $k(\cdot)$ at stage $t$.
\end{itemize}
Before computing a new sample at stage $t+1$, each weight $p_o$ is updated according to the formula $p_o\cdot (h^{(t)})^{-m^{{(t+1)}}}\cdot k\big(\frac{\widehat{a}_{i}^{(t)}-a_{i,o}^{(t)}}{h^{(t)}}\big)$, and then normalized. Further, a random number $\alpha$ is drawn anew.\\
Using this procedure, the conditional density $g_{i}^{(t+1)}$ of the accumulation rate of bin $i$ at stage $t+1$, given $\widehat{a}_{i}^{(1)},\ldots,\widehat{a}_{i}^{(t)}$, can be estimated by
\begin{linenomath}
\begin{equation*}
	\widehat{g}_i^{(t+1)} \big(\widehat{a}_{i}^{(t+1)} \big| \widehat{a}_{i}^{(1)},\ldots,\widehat{a}_{i}^{(t)}\big) = \sum_{o=1}^{N_o} p_o \cdot (h^{(N_o)})^{-m^{(t+1)}}\cdot k\bigg(\frac{\widehat{a}_{i}^{(t+1)}-a_{i,o}^{(t+1)}}{h^{N_o}}\bigg).
\end{equation*}
\end{linenomath} 
Within this approach, every new trajectory starts at $\widehat{a}_{i}^{(1)}$, and new samples $\widehat{a}_{i}^{(t+1)}$ are generated according to the density $\widehat{g}_i^{(t+1)}$, for $t=1,\ldots,T-1$. At the end of the procedure at stage $T$, a new trajectory $(\widehat{a}_{i}^{(1)},\ldots,\widehat{a}_{i}^{(T)})$ has been generated from the initial data.

\begin{figure}[h!]
\centering
\begin{subfigure}[b]{0.49\textwidth}
	\centering
	\includegraphics[width=\textwidth]{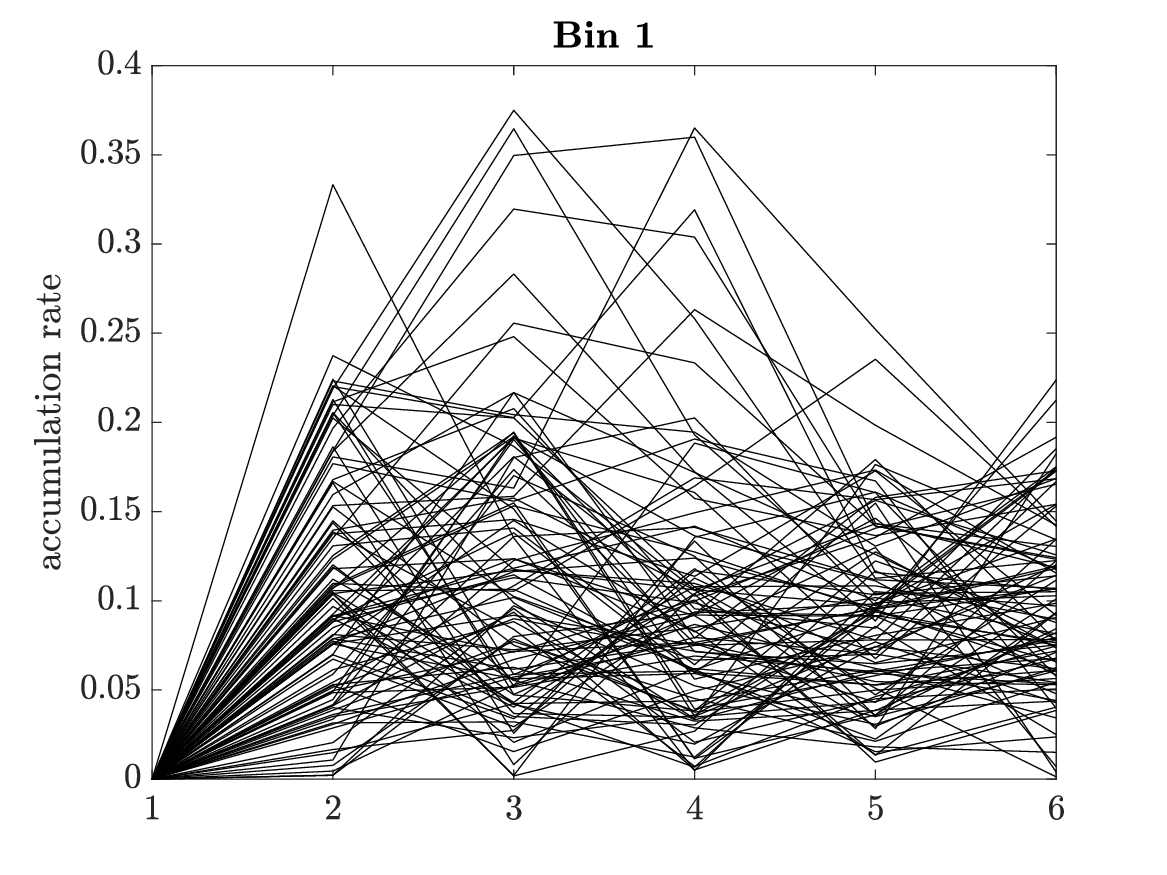}
\end{subfigure} \hfill
\begin{subfigure}[b]{0.49\textwidth}
	\centering
	\includegraphics[width=\textwidth]{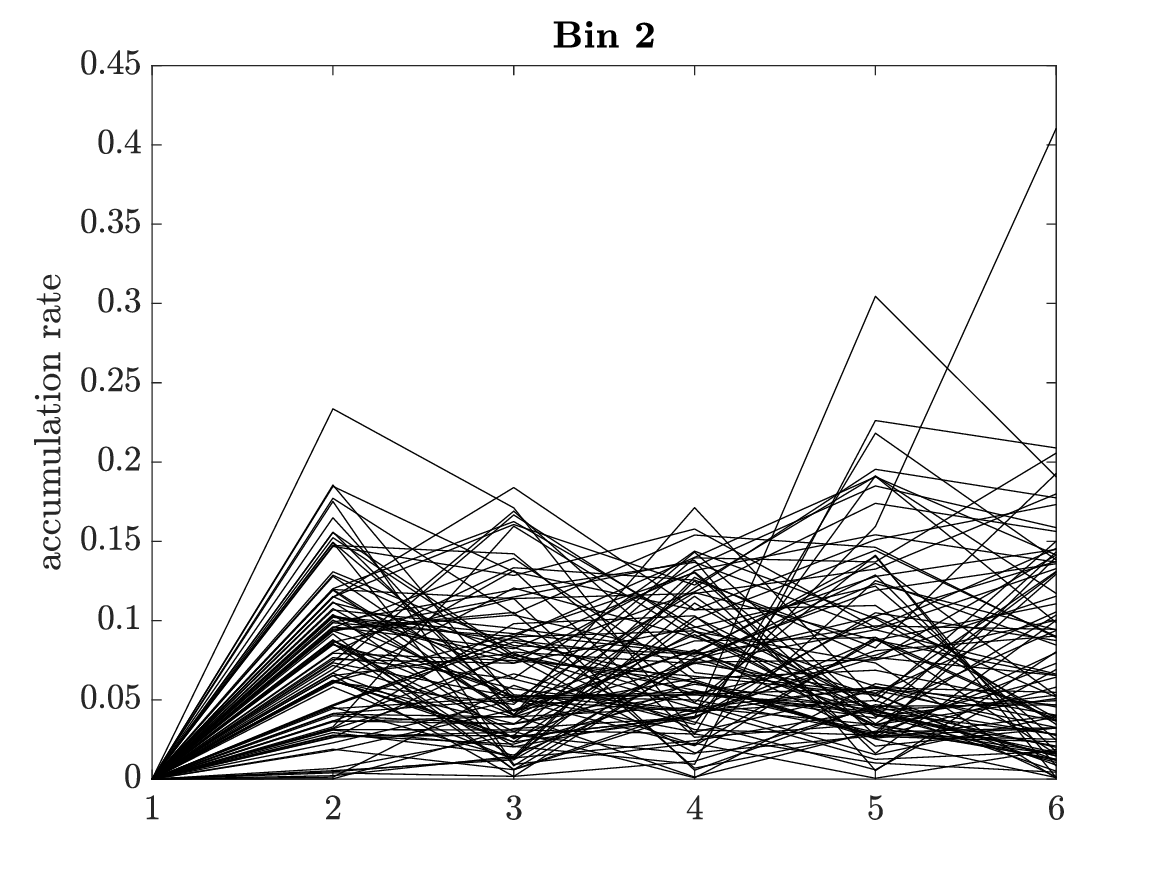}
\end{subfigure}
\begin{subfigure}[b]{0.49\textwidth}
	\centering
	\includegraphics[width=\textwidth]{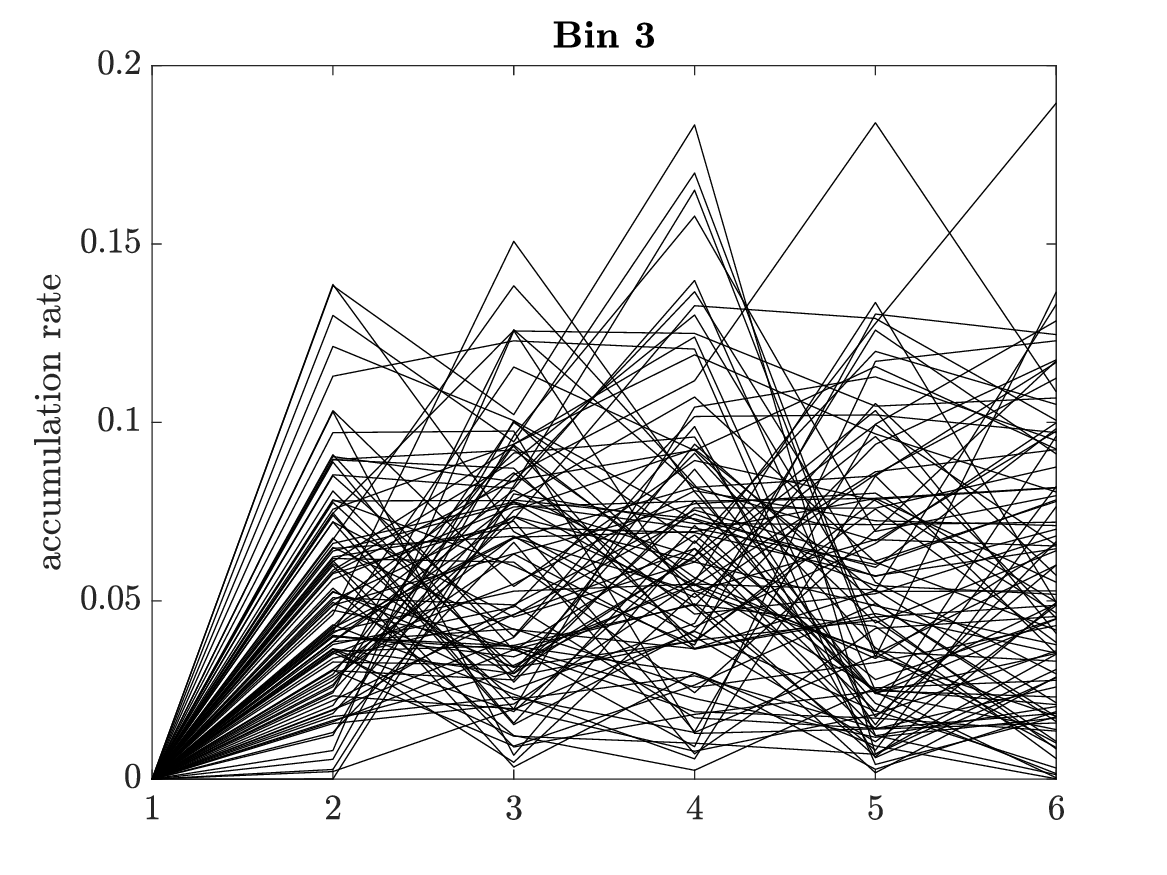}
\end{subfigure} \hfill
\begin{subfigure}[b]{0.49\textwidth}
	\centering
	\includegraphics[width=\textwidth]{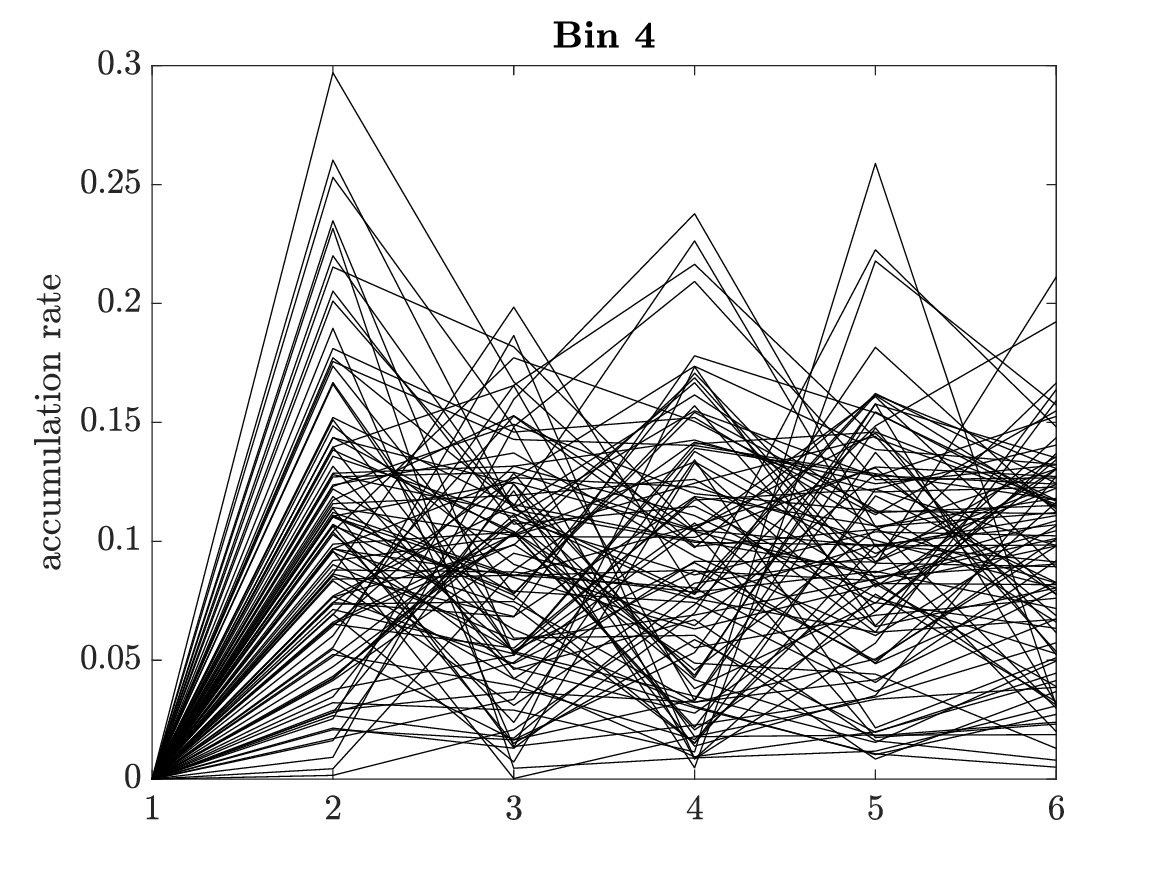}
\end{subfigure}
\begin{subfigure}[b]{0.49\textwidth}
	\centering
	\includegraphics[width=\textwidth]{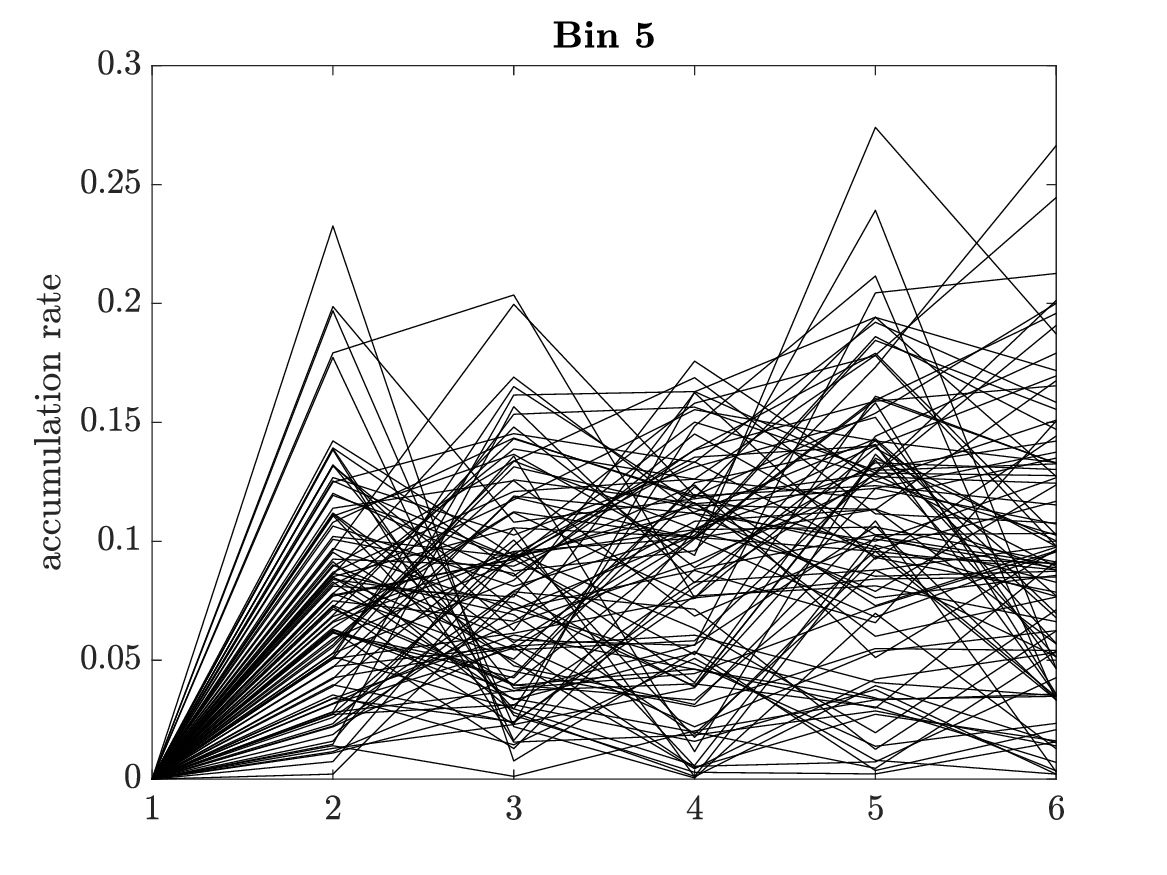}
\end{subfigure} \hfill
\begin{subfigure}[b]{0.49\textwidth}
	\centering
	\includegraphics[width=\textwidth]{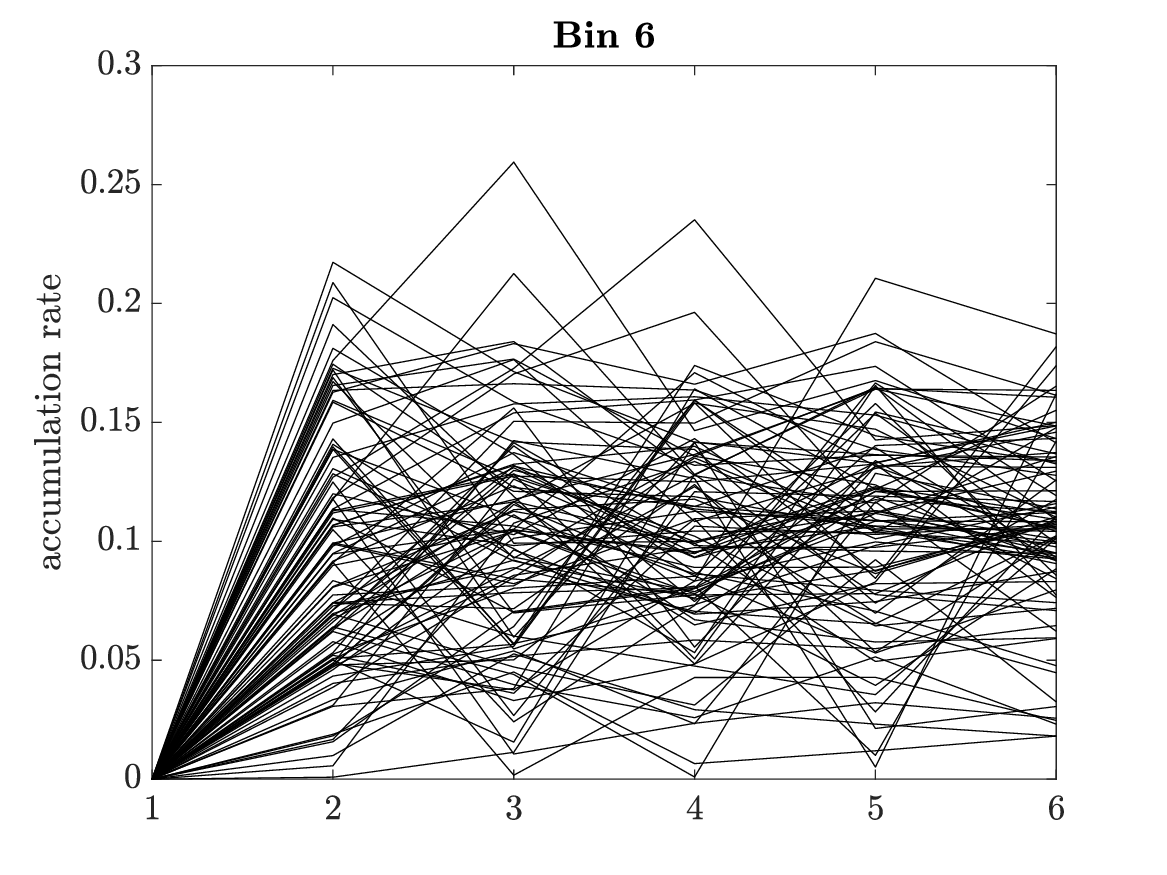}
\end{subfigure}
\caption{For each of the six considered bins, one hundred trajectories on the accumulation rate of waste generated from historical data through the conditional density estimation process are depicted. The stages are represented on the horizontal axis.}
\label{fig_trajectories}
\end{figure}

We set $a_{i,o}^{(1)}=0=\widehat{a}_{i}^{(1)}$ for all $i=1,\ldots,N$, $o=1,\ldots,N_o$ because no increase of waste at the first stage of the time horizon is assumed, and $m^{(t)}=N$ for all $t=1,\ldots,T$ since, at each node, the dimension of the state corresponds to the total number of bins. Furthermore, as suggested in \cite{KirPflPic2020}, the kernel $k(\cdot)$ is set to be logistic. Figure \ref{fig_trajectories} shows one hundred trajectories of the accumulation rate in six different bins, generated according to the conditional density estimation process described so far.

Secondly, we apply a dynamic stochastic approximation algorithm to generate a candidate scenario tree (see \cite{PflPic2016} for details). Starting from an initial guess of a tree with a prescribed branching structure, at every iteration of the procedure a new sample path is generated according to the conditional density estimation process discussed above. The algorithm finds one possible sequence of nodes in the scenario tree whose distance between the states of those nodes and the generated sample is minimal. Thus, the states of those nodes are updated with the values of the generated sample and the others remain unchanged. Then, the algorithm calculates the conditional probabilities to reach each node of the tree starting from its root, and it stops when all the iterations, whose number is decided in advance, have been performed.

The scenario tree generation procedure described so far has been implemented in Julia, relying on the package ScenTrees.jl (see \cite{KirPicPfl2020JULIA}). The number of iterations for the stochastic approximation process has been set to 10000.

\clearpage
\emph{C)} \emph{In-sample stability}\\
In this section, we carry out an in-sample stability analysis (see \cite{KauWal2007}).

\begin{figure}[ht!]
\centering
\includegraphics[width=0.7\textwidth]{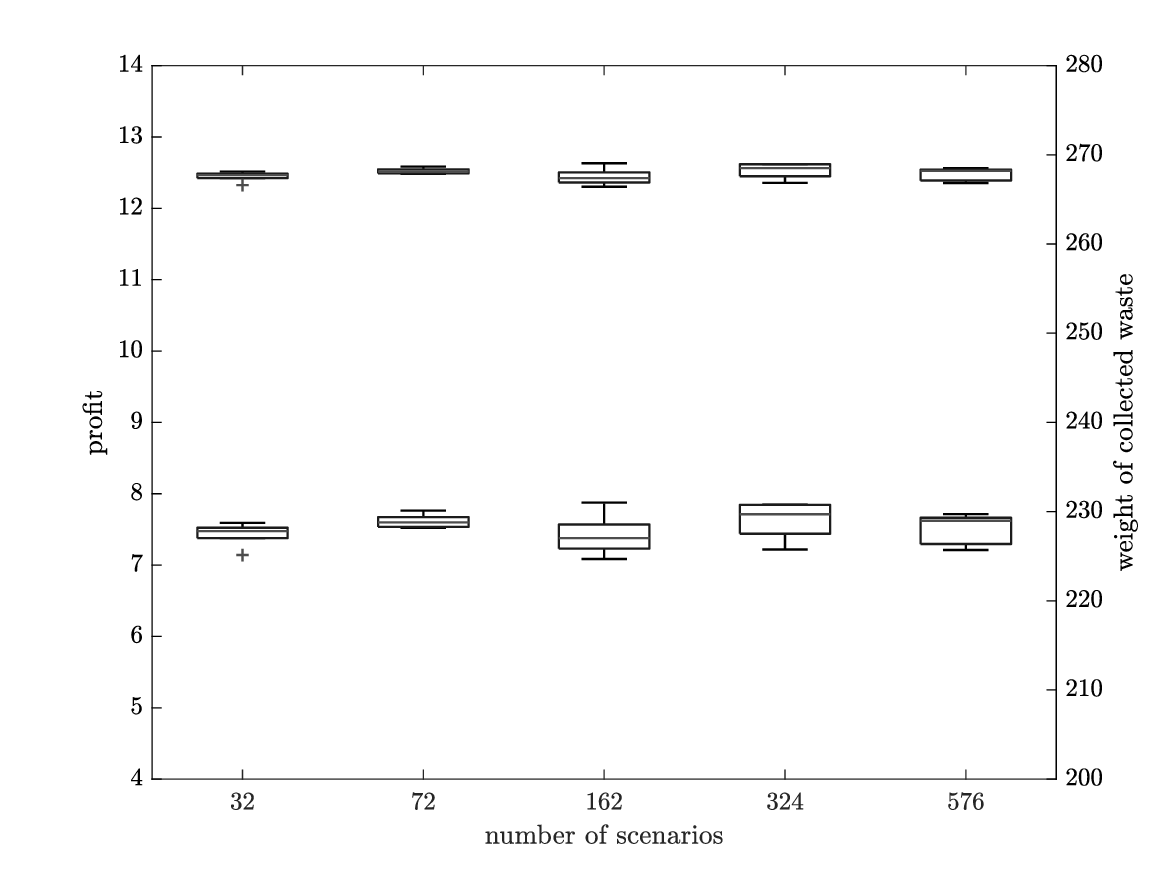}
\caption{Box-plots of objective function value (below, left-hand scale) and of weight of collected waste (above, right-hand scale) over 5 runs of scenario trees with increasing cardinality.}
\label{fig_insample_stability}
\end{figure}

In Table \ref{tab_insample_stability} we report average results obtained by solving model $\mathcal{M}$ over five runs on \textit{inst\_9\_1}, with increasing size of the scenario tree. Box-plots of objective function and of weight of collected waste are depicted in Figure \ref{fig_insample_stability}.

\begin{table}[h!]
\centering
\resizebox{\textwidth}{!}{
	\begin{tabular}{ccccccc}\toprule
		Scenarios & Branching structure & Profit (\euro) & Weight of waste (kg) & Distance (km) & CPU time (s) & Multistage distance \\ \hline
		32 & [1 2 2 2 2 2] & 7.43 & 267.58 & 72.84 & 55.22 & 0.063\\
		72 & [1 3 3 2 2 2] & 7.61 & 268.18 & 72.84 & 166.24 & 0.046\\
		162 & [1 3 3 3 3 2] & 7.42 & 267.52 & 72.84 & 1165.00 & 0.034\\
		324 & [1 4 3 3 3 3] & 7.63 & 268.22 & 72.84 & 7109.18 & 0.027\\
		576 & [1 4 4 4 3 3] & 7.50 & 267.81 & 72.84 & 51928.64 & 0.020\\
		1024 & [1 4 4 4 4 4] & \multicolumn{4}{c}{Not solved to optimality within 24 hours} & 0.016\\
		\bottomrule
\end{tabular}}
\caption{Average results on the in-sample stability analysis over five runs on scenario trees with increasing size. The results are drawn from model $\mathcal{M}$ on \textit{inst\_9\_1}.} \label{tab_insample_stability}
\end{table}

Since various indicators (profit, weight of collected waste, total travelled distance) do not vary significantly when increasing the size of the tree, we conclude that the methodology we applied to generate scenario trees is stable even with small trees. Besides, the multistage distance (see the last column of Table \ref{tab_insample_stability}), is throughout close to zero, due to the minimization of the distance in the dynamic stochastic approximation algorithm (see \cite{PflPic2016}. On the other hand, the computational time increases considerably, when increasing the size of the tree.

\begin{figure}[h!]
\centering
\begin{subfigure}[b]{0.49\textwidth}
	\centering
	\includegraphics[width=\textwidth]{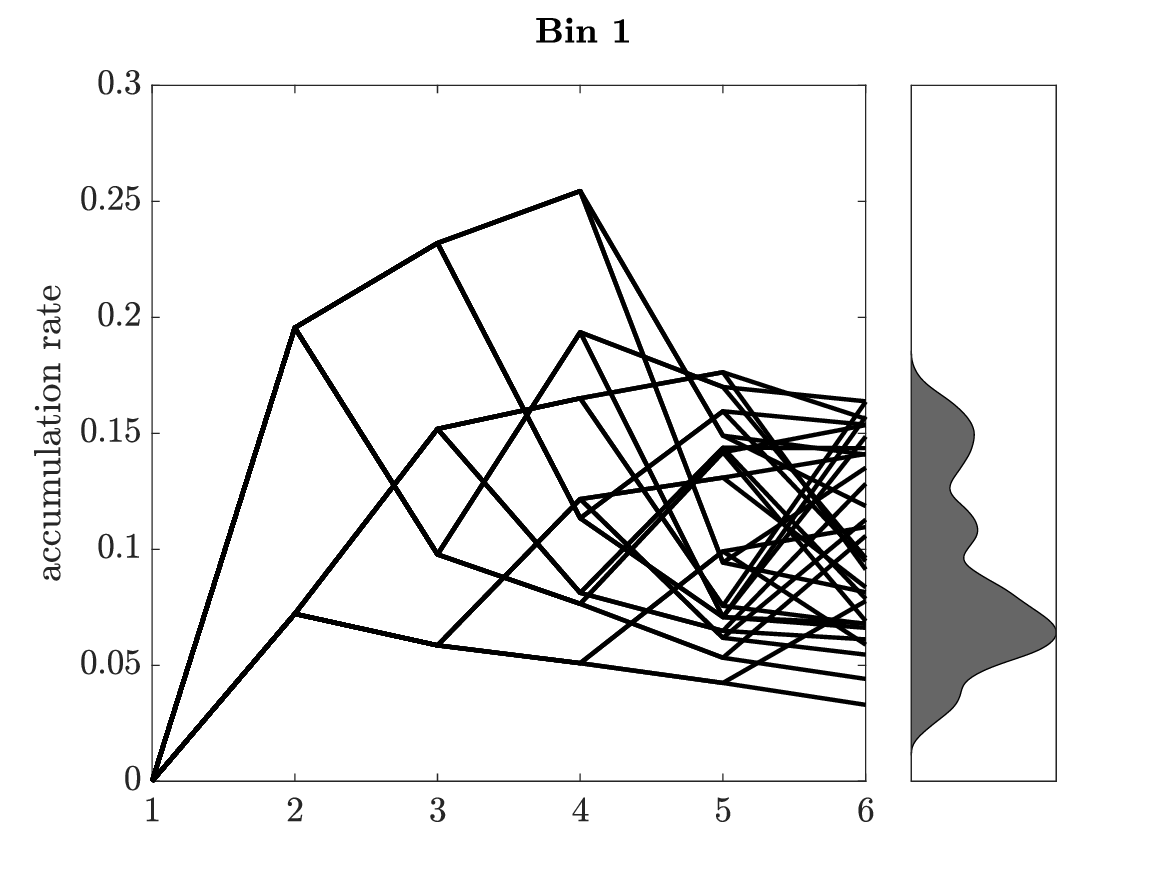}
\end{subfigure} \hfill
\begin{subfigure}[b]{0.49\textwidth}
	\centering
	\includegraphics[width=\textwidth]{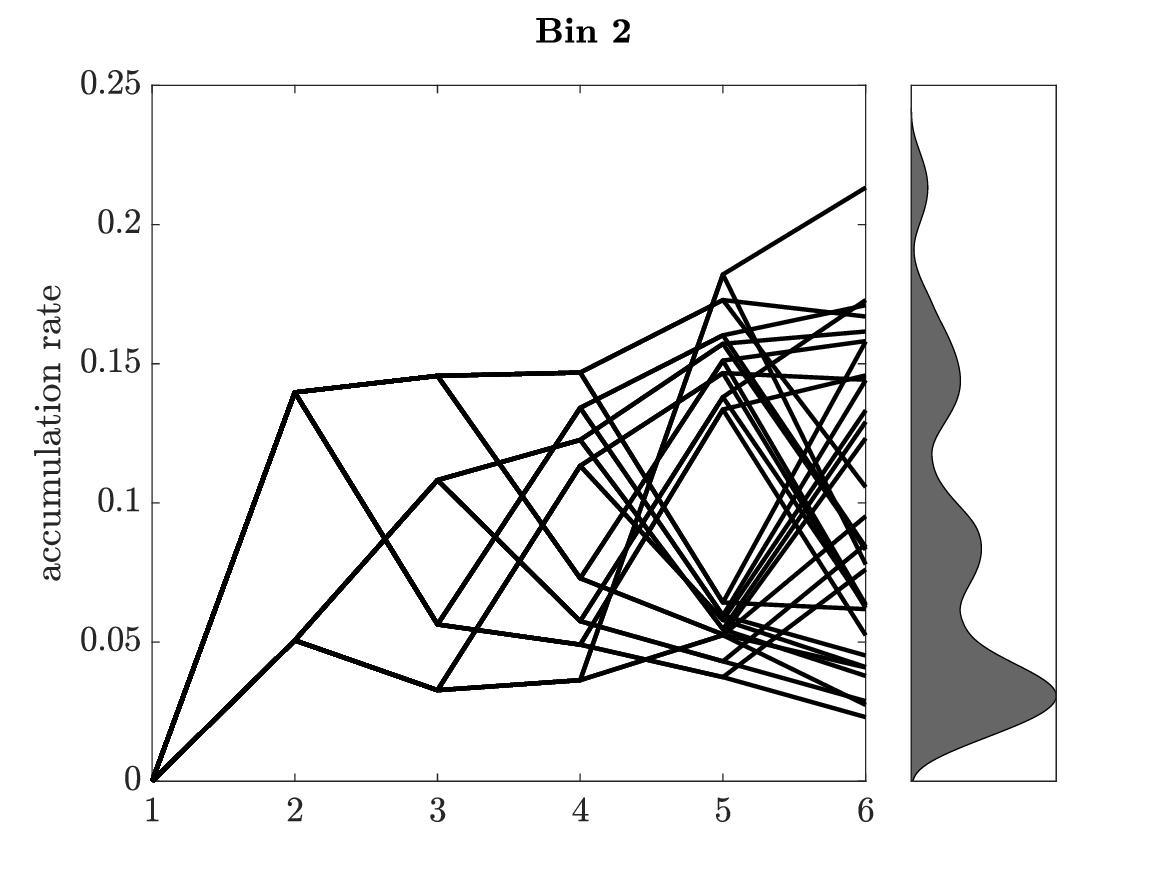}
\end{subfigure}
\begin{subfigure}[b]{0.49\textwidth}
	\centering
	\includegraphics[width=\textwidth]{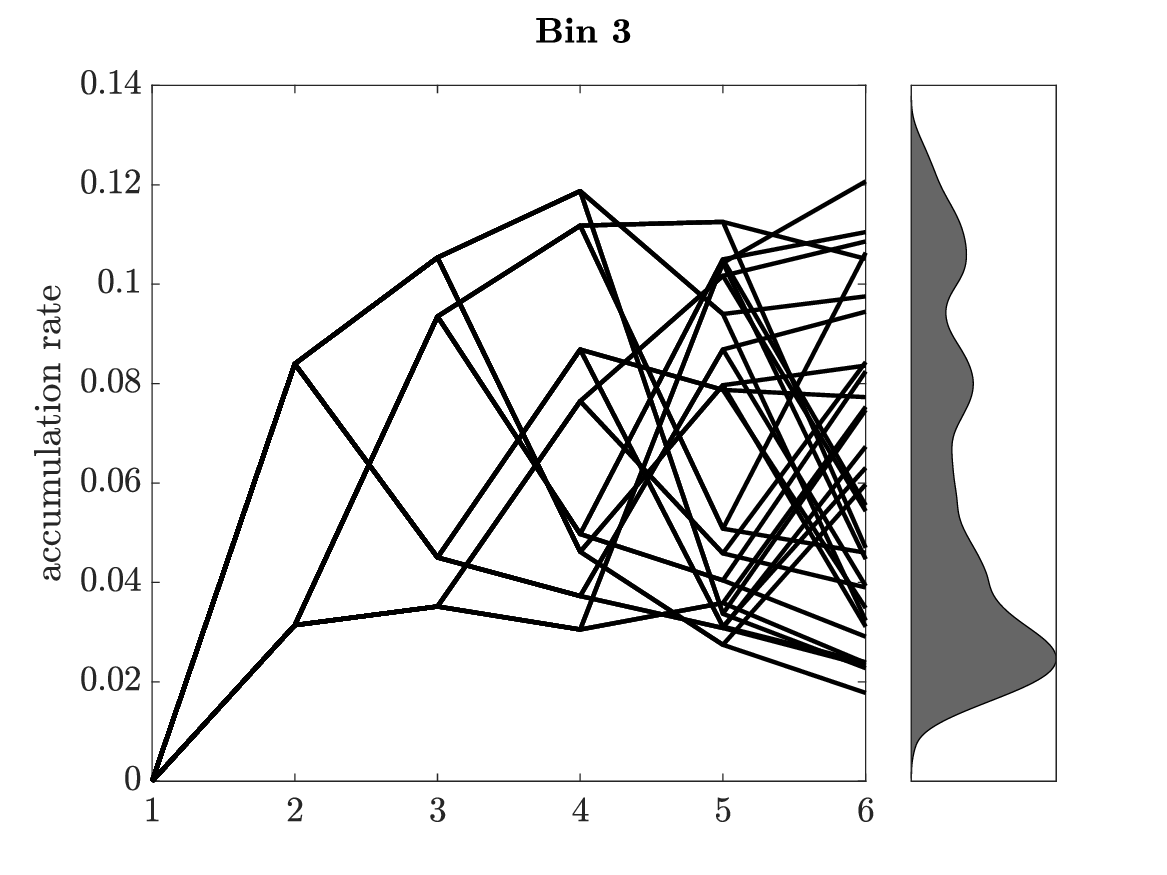}
\end{subfigure} \hfill
\begin{subfigure}[b]{0.49\textwidth}
	\centering
	\includegraphics[width=\textwidth]{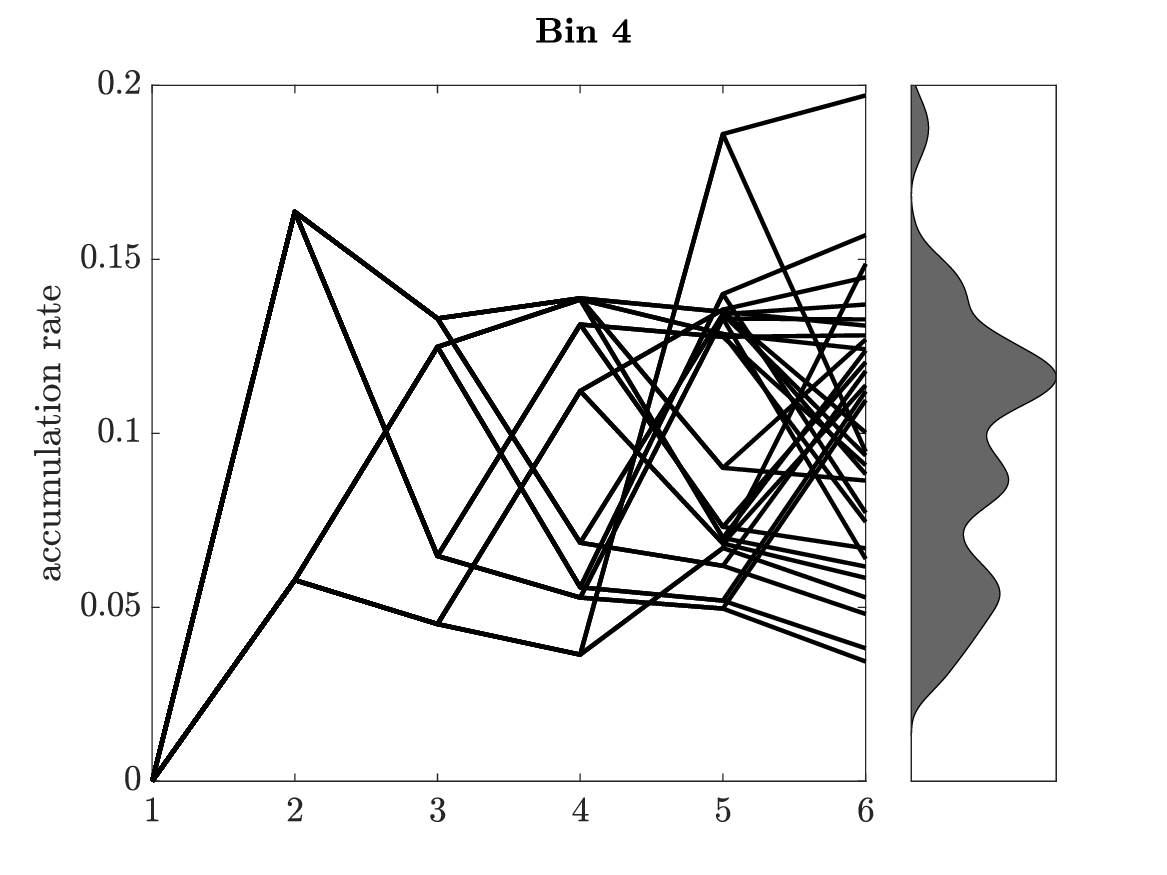}
\end{subfigure}
\begin{subfigure}[b]{0.49\textwidth}
	\centering
	\includegraphics[width=\textwidth]{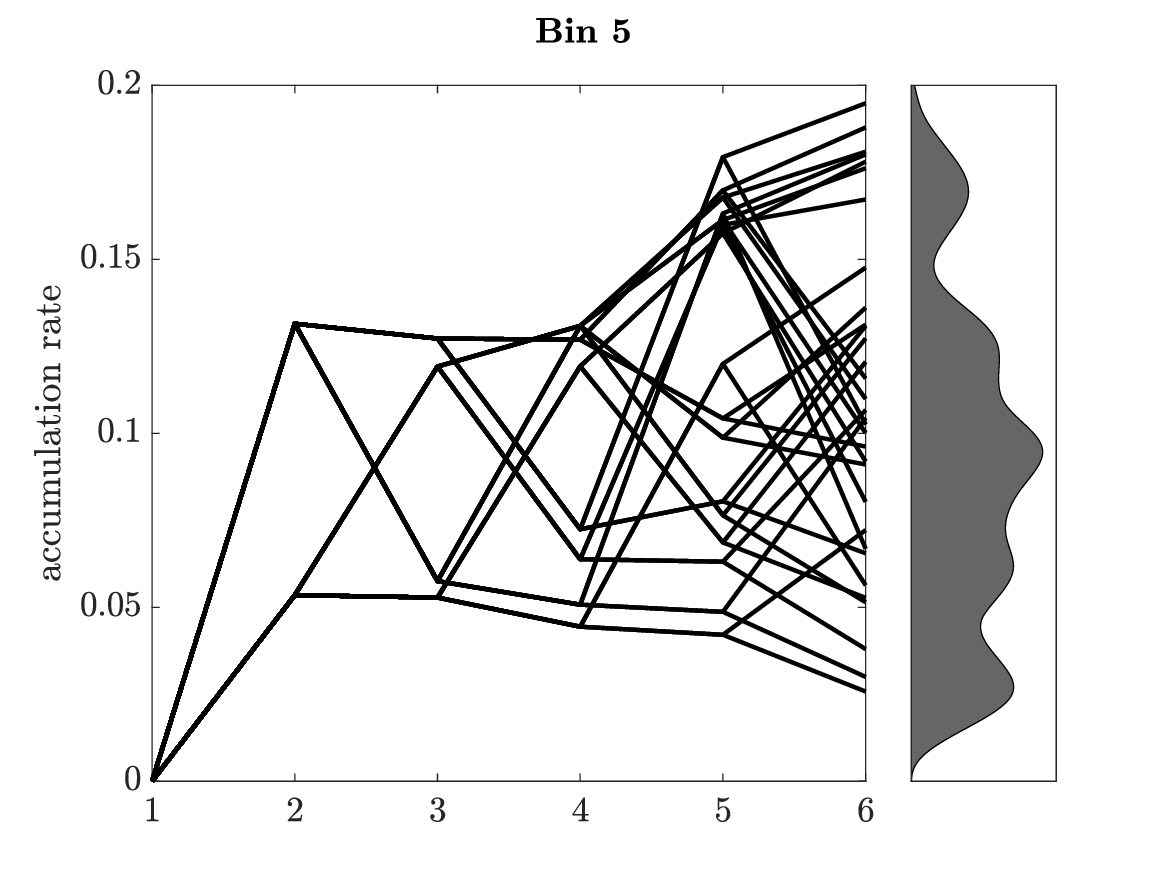}
\end{subfigure} \hfill
\begin{subfigure}[b]{0.49\textwidth}
	\centering
	\includegraphics[width=\textwidth]{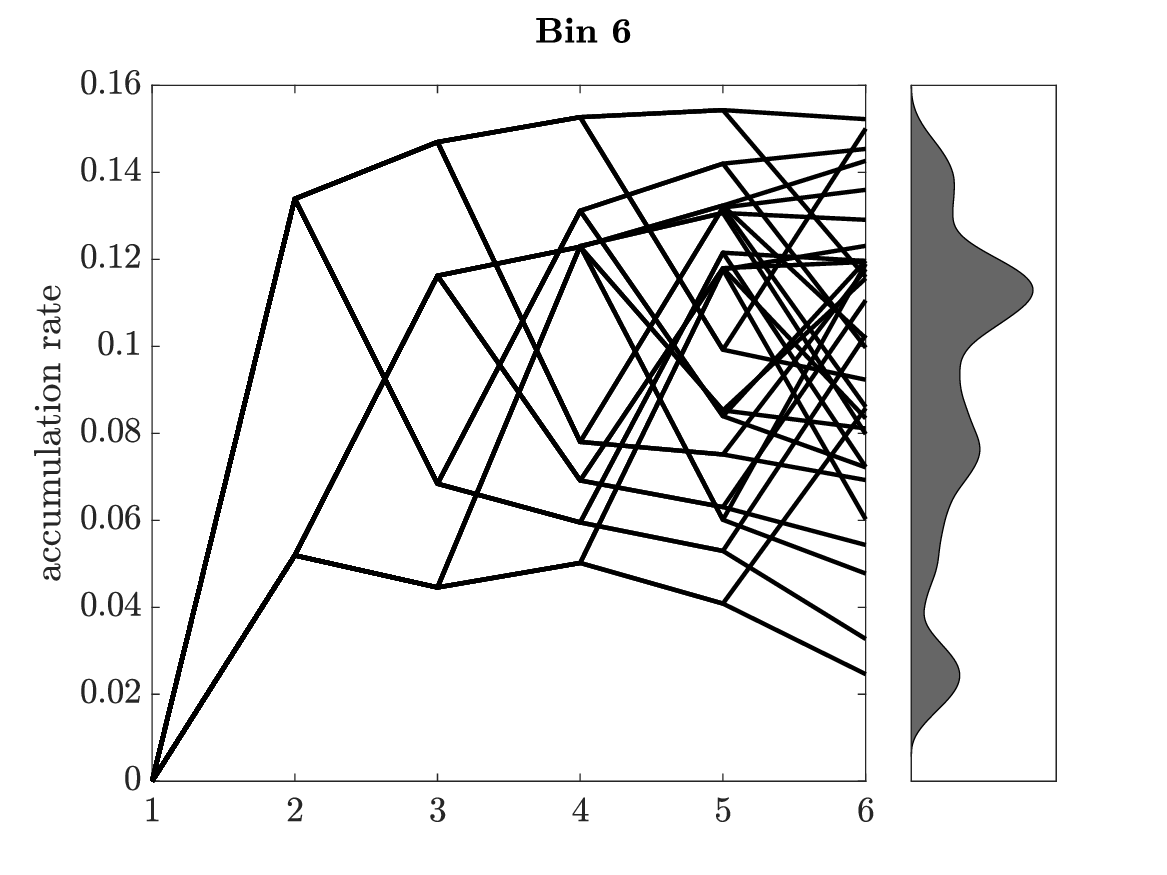}
\end{subfigure}
\caption{Examples of six-stage scenario trees of the accumulation rate of waste in the same six bins as in Figure \ref{fig_trajectories}. The corresponding probability distribution is depicted on the right of each plot.}
\label{fig_scenario_trees_four}
\end{figure}

For all of these reasons, we decide to consider a scenario tree of size $S=32$, with 63 nodes. In Figure \ref{fig_scenario_trees_four} we depict six binary scenario trees of six different bins with the corresponding probability distributions generated from the dynamic stochastic approximation algorithm.

\clearpage
\emph{D)} \emph{Stochastic measures (detailed results for small instances)}

\begin{table}[h!]
\centering
\resizebox{0.9\textwidth}{!}{
	\begin{tabular}{lllllllllll}\toprule
		& \textit{inst\_1\_9} & \textit{inst\_2\_9} & \textit{inst\_3\_9} & \textit{inst\_4\_9} & \textit{inst\_5\_9} & \textit{inst\_6\_9} & \textit{inst\_7\_9} & \textit{inst\_8\_9} & \textit{inst\_9\_9} & \textit{inst\_10\_9} \\ \hline
		$RP$ & 9.36 & 11.92 & 31.58 & 32.66 & 2.48 & 4.27 & 30.76 & 22.90 & 2.72 &	32.96\\
		$EV$ & 23.79	& 14.43 & 38.23 & 38.58 & 17.24 & 18.10 & 36.85 & 45.80 & 18.60 & 42.82\\
		$WS$ & 17.63 & 17.46 & 45.50 & 37.94 & 16.40 & 15.88 & 35.98 & 25.40 & 24.04 & 45.68\\ \hline
		$\%EVPI$ & 88\%	& 46\%	& 44\%	& 16\% &	562\% & 272\% & 17\% & 11\% & 783\%	& 39\%\\ \hline
		$\%VSS^1$ & $\infty$ & $\infty$ & $\infty$ & 77\% & $\infty$ & $\infty$ & $\infty$ & $\infty$ & $\infty$ & $\infty$\\
		$\%VSS^2$ & $\infty$ & $\infty$ & $\infty$ & 77\% & $\infty$ & $\infty$ & $\infty$ & $\infty$ & $\infty$ & $\infty$\\
		$\%VSS^3$ & $\infty$ & $\infty$ & $\infty$ & $\infty$ & $\infty$ & $\infty$  & $\infty$ & $\infty$ & $\infty$ & $\infty$\\
		$\%VSS^4$ & $\infty$ & $\infty$ & $\infty$ & $\infty$ & $\infty$ & $\infty$  & $\infty$ & $\infty$ & $\infty$ & $\infty$\\
		$\%VSS^5$ & $\infty$ & $\infty$ & $\infty$ & $\infty$ & $\infty$ & $\infty$  & $\infty$ & $\infty$ & $\infty$ & $\infty$\\ \hline
		$\%MLUSS^1$ & $\infty$ & $\infty$ & $\infty$ & 77\% & $\infty$ & $\infty$ & $\infty$ & $\infty$ & $\infty$ & $\infty$\\
		$\%MLUSS^2$ & $\infty$ & $\infty$ & $\infty$ & 77\% & $\infty$ & $\infty$ & $\infty$ & $\infty$ & $\infty$ & $\infty$\\
		$\%MLUSS^3$ & $\infty$ & $\infty$ & $\infty$ & $\infty$ & $\infty$ & $\infty$  & $\infty$ & $\infty$ & $\infty$ & $\infty$\\
		$\%MLUSS^4$ & $\infty$ & $\infty$ & $\infty$ & $\infty$ & $\infty$ & $\infty$  & $\infty$ & $\infty$ & $\infty$ & $\infty$\\
		$\%MLUSS^5$ & $\infty$ & $\infty$ & $\infty$ & $\infty$ & $\infty$ & $\infty$  & $\infty$ & $\infty$ & $\infty$ & $\infty$\\ \hline
		$\%MLUDS^1$ & 0\%	& 0\%	& 0\% &	77\% &	0\%	& 0\%	& 0\%	& 0\%	& 0\%	& 0\%\\
		$\%MLUDS^2$ & 546\% &	500\%	& 179\%	& 77\%	& 1992\% &	1131\%	& 174\%	 & 235\%	 & 1760\%	& 150\%\\
		$\%MLUDS^3$ & 546\% &	500\%	& 179\%	& 77\% &	1992\%	&1131\%	 & 174\% & 235\%	& 1760\% &	150\%\\
		$\%MLUDS^4$ & 546\% &	500\%	& 179\%	& 77\% &	1992\%	&1131\%	 & 174\% & 235\%	& 1760\% &	150\%\\
		$\%MLUDS^5$ & 546\% &	500\%	& 179\%	& 147\% &	1992\%	&1131\%	 & 174\% & 235\%	& 1760\% &	150\%\\
		\bottomrule
\end{tabular}}
\caption{Detailed results of $RP$, $EV$, $WS$ and of stochastic measures $\%EVPI$, $\%VSS^t$, $\%MLUSS^t$, $\%MLUDS^t$, for $1\leq t \leq 5$. The values in percentage denote the gap with respect to the corresponding $RP$ problem. The results refer to the instances with 9 bins.} \label{tab_res_9bins}
\end{table}

\begin{table}[h!]
\centering
\resizebox{0.9\textwidth}{!}{
	\begin{tabular}{lllllllllll}\toprule
		& \textit{inst\_1\_10} & \textit{inst\_2\_10} & \textit{inst\_3\_10} & \textit{inst\_4\_10} & \textit{inst\_5\_10} & \textit{inst\_6\_10} & \textit{inst\_7\_10} & \textit{inst\_8\_10} & \textit{inst\_9\_10} & \textit{inst\_10\_10} \\ \hline
		$RP$ & 14.57 & 	25.97	& 53.88	& 54.09	& 16.41	& 32.07	& 41.28	& 35.71	& 33.16	& 36.59\\
		$EV$ & 32.06 &	28.48 &	53.88	& 58.82	& 22.50	& 32.07	& 48.12	& 40.65	& 40.45	& 36.59\\
		$WS$ & 34.42 & 32.40	& 57.08	& 63.13	& 29.61	& 35.38	& 47.91	& 49.40	& 43.67	& 38.10\\
		\hline
		$\%EVPI$ & 136\% &	25\%	& 6\%	 & 17\% &	80\%	 & 10\% &	16\% &	38\%	 & 32\%	 & 4\% \\
		\hline
		$\%VSS^1$ & $\infty$ & $\infty$ & $\infty$ & $\infty$ & $\infty$ & $\infty$ & 0\% & $\infty$ & $\infty$ & $\infty$\\
		$\%VSS^2$ & $\infty$ & $\infty$ & $\infty$ & $\infty$ & $\infty$ & $\infty$ & 55\% & $\infty$ & $\infty$ & $\infty$\\
		$\%VSS^3$ & $\infty$ & $\infty$ & $\infty$ & $\infty$ & $\infty$ & $\infty$  & 55\% & $\infty$ & $\infty$ & $\infty$\\
		$\%VSS^4$ & $\infty$ & $\infty$ & $\infty$ & $\infty$ & $\infty$ & $\infty$  & $\infty$ & $\infty$ & $\infty$ & $\infty$\\
		$\%VSS^5$ & $\infty$ & $\infty$ & $\infty$ & $\infty$ & $\infty$ & $\infty$  & $\infty$ & $\infty$ & $\infty$ & $\infty$\\
		\hline
		$\%MLUSS^1$ & $\infty$ & $\infty$ & $\infty$ & $\infty$ & $\infty$ & $\infty$ & 0\% & $\infty$ & $\infty$ & $\infty$\\
		$\%MLUSS^2$ & $\infty$ & $\infty$ & $\infty$ & $\infty$ & $\infty$ & $\infty$ & 55\% & $\infty$ & $\infty$ & $\infty$\\
		$\%MLUSS^3$ & $\infty$ & $\infty$ & $\infty$ & $\infty$ & $\infty$ & $\infty$  & 55\% & $\infty$ & $\infty$ & $\infty$\\
		$\%MLUSS^4$ & $\infty$ & $\infty$ & $\infty$ & $\infty$ & $\infty$ & $\infty$  & $\infty$ & $\infty$ & $\infty$ & $\infty$\\
		$\%MLUSS^5$ & $\infty$ & $\infty$ & $\infty$ & $\infty$ & $\infty$ & $\infty$  & $\infty$ & $\infty$ & $\infty$ & $\infty$\\ \hline
		$\%MLUDS^1$ & 0\% & 	0\%	& 0\%	& 0\%	& 0\%	 &0\%&	0\% &	0\%	& 0\%	& 0\%\\
		$\%MLUDS^2$ & 316\%	& 229\%	&	95\%	&	87\%	&	348\%	&	160\%	&	55\%	&	153\%	&	157\%	&	147\%\\
		$\%MLUDS^3$ & 316\%	& 229\%	&	95\%	&	87\%	&	348\%	&	160\%	&	55\%	&	153\%	&	157\%	&	147\%\\
		$\%MLUDS^4$ & 316\%	& 229\%	&	95\%	&	87\%	&	348\%	&	160\%	&	55\%	&	153\%	&	157\%	&	147\%\\
		$\%MLUDS^5$ & 316\%	& 229\%	&	95\%	&	87\%	&	348\%	&	160\%	&	55\%	&	153\%	&	157\%	&	147\%\\
		\bottomrule
\end{tabular}}
\caption{Detailed results of $RP$, $EV$, $WS$ and of stochastic measures $\%EVPI$, $\%VSS^t$, $\%MLUSS^t$, $\%MLUDS^t$, for $1\leq t \leq 5$. The values in percentage denote the gap with respect to the corresponding $RP$ problem. The results refer to the instances with 10 bins.} \label{tab_res_10bins}
\end{table}

\begin{table}[h!]
\centering
\resizebox{0.9\textwidth}{!}{
	\begin{tabular}{lllllllllll}\toprule
		& \textit{inst\_1\_11} & \textit{inst\_2\_11} & \textit{inst\_3\_11} & \textit{inst\_4\_11} & \textit{inst\_5\_11} & \textit{inst\_6\_11} & \textit{inst\_7\_11} & \textit{inst\_8\_11} & \textit{inst\_9\_11} & \textit{inst\_10\_11} \\ \hline
		$RP$ & 30.83 &	38.46 &	64.24	 &33.12	&46.72	&50.21	&60.99	&29.87	&15.73	&42.31\\
		$EV$ &32.38	&41.03	&66.18	&33.12	&49.84	&53.15	&61.94	&34.10	&26.57	&52.39\\
		$WS$ & 40.32 &	41.96	& 65.57	& 46.96	& 51.60	& 52.02	& 62.36	& 36.42	& 25.93	& 48.32\\
		\hline
		$\%EVPI$ & 31\%	& 9\%	& 2\%	& 42\%	& 10\%	& 4\%	& 2\%	& 22\%	& 65\%	& 14\%\\
		\hline
		$\%VSS^1$ & $\infty$ & $\infty$ & $\infty$ & $\infty$ & $\infty$ & $\infty$ & $\infty$ & $\infty$ & $\infty$ & $\infty$\\
		$\%VSS^2$ & $\infty$ & $\infty$ & $\infty$ & $\infty$ & $\infty$ & $\infty$ & $\infty$ & $\infty$ & $\infty$ & $\infty$\\
		$\%VSS^3$ & $\infty$ & $\infty$ & $\infty$ & $\infty$ & $\infty$ & $\infty$  & $\infty$ & $\infty$ & $\infty$ & $\infty$\\
		$\%VSS^4$ & $\infty$ & $\infty$ & $\infty$ & $\infty$ & $\infty$ & $\infty$  & $\infty$ & $\infty$ & $\infty$ & $\infty$\\
		$\%VSS^5$ & $\infty$ & $\infty$ & $\infty$ & $\infty$ & $\infty$ & $\infty$  & $\infty$ & $\infty$ & $\infty$ & $\infty$\\
		\hline
		$\%MLUSS^1$ & $\infty$ & $\infty$ & $\infty$ & $\infty$ & $\infty$ & $\infty$ & $\infty$ & $\infty$ & $\infty$ & $\infty$\\
		$\%MLUSS^2$ & $\infty$ & $\infty$ & $\infty$ & $\infty$ & $\infty$ & $\infty$ & $\infty$ & $\infty$ & $\infty$ & $\infty$\\
		$\%MLUSS^3$ & $\infty$ & $\infty$ & $\infty$ & $\infty$ & $\infty$ & $\infty$  & $\infty$ & $\infty$ & $\infty$ & $\infty$\\
		$\%MLUSS^4$ & $\infty$ & $\infty$ & $\infty$ & $\infty$ & $\infty$ & $\infty$  & $\infty$ & $\infty$ & $\infty$ & $\infty$\\
		$\%MLUSS^5$ & $\infty$ & $\infty$ & $\infty$ & $\infty$ & $\infty$ & $\infty$  & $\infty$ & $\infty$ & $\infty$ & $\infty$\\ \hline
		$\%MLUDS^1$ & 0\% & 	0\%	& 0\%	& 0\%	& 0\%	 &0\%&	0\% &	0\%	& 0\%	& 0\%\\
		$\%MLUDS^2$ & 200\%	& 155\%	&	100\%	&	156\%	&	129\%	&	124\%	&	93\%	&	166\%	&	326\%	&	131\%\\
		$\%MLUDS^3$ & 200\%	& 155\%	&	100\%	&	156\%	&	129\%	&	124\%	&	93\%	&	166\%	&	326\%	&	131\%\\
		$\%MLUDS^4$ & 200\%	& 155\%	&	100\%	&	156\%	&	129\%	&	124\%	&	93\%	&	166\%	&	326\%	&	131\%\\
		$\%MLUDS^5$ & 200\%	& 155\%	&	100\%	&	156\%	&	129\%	&	124\%	&	93\%	&	166\%	&	326\%	&	131\%\\
		\bottomrule
\end{tabular}}
\caption{Detailed results of $RP$, $EV$, $WS$ and of stochastic measures $\%EVPI$, $\%VSS^t$, $\%MLUSS^t$, $\%MLUDS^t$, for $1\leq t \leq 5$. The values in percentage denote the gap with respect to the corresponding $RP$ problem. The results refer to the instances with 11 bins.} \label{tab_res_11bins}
\end{table}

\clearpage
\emph{E)} \emph{Performance of the rolling horizon approach (detailed results for small instances)}

\begin{table}[h!]
\centering
\resizebox{0.9\textwidth}{!}{
	\begin{tabular}{lllllllllll}\toprule
		& \textit{inst\_1\_9} & \textit{inst\_2\_9} & \textit{inst\_3\_9} & \textit{inst\_4\_9} & \textit{inst\_5\_9} & \textit{inst\_6\_9} & \textit{inst\_7\_9} & \textit{inst\_8\_9} & \textit{inst\_9\_9} & \textit{inst\_10\_9} \\ \hline
		$W$ & \multicolumn{10}{c}{\bf Profit reduction ($\%$)}\\
		1 & 11\% & 	$\infty$ &	37\% &	8\% &	74\% &	0\%	& 32\%	& $\infty$	& 0\%	& 7\%\\
		2 & 11\% & 	95\% &	37\% &	8\% &	74\%	& 0\%	& 32\% &	54\%	& 0\%	& 7\%\\
		3 & 11\% &	95\% &	37\% &	8\%	& 74\%	& 0\%	& 32\%	& 54\%	& 0\%	& 7\%\\
		4 & 11\%	 & 0\%	& 37\%	& 8\%	& 74\%	& 0\%	& 32\%	& 0\%	& 0\%	& 7\%\\
		\hline
		$W$ & \multicolumn{10}{c}{\bf Computational time reduction ($\%$)}\\
		1 & 94\% & 	99\%	& 99\%	& 94\%	& 99\%	& 99\%	& 100\%	& 95\%	& 100\%	& 100\%\\
		2 & 85\%	 & 76\%	& 96\%	& 81\%	& 98\%	& 98\%	& 99\%	& 79\%	& 100\%	& 97\%\\
		3 & 67\%	 & 49\%	& 90\%	& 49\%	& 83\%	& 96\%	& 96\%	& 49\%	& 99\%	& 91\%\\
		4 & 40\% & $-9\%$ & 	75\%	& 5\%	& 85\%	& 79\%	& 94\%	& $-36\%$	& 98\%	& 67\%\\
		\bottomrule
\end{tabular}}
\caption{Detailed results on the performance of the rolling horizon approach, in terms of reduction of the profit and of the CPU time when compared to the $RP$ problem. The results refer to the instances with 9 bins.} \label{tab_RH_9bins}
\end{table}

\begin{table}[h!]
\centering
\resizebox{0.9\textwidth}{!}{
	\begin{tabular}{lllllllllll}\toprule
		& \textit{inst\_1\_10} & \textit{inst\_2\_10} & \textit{inst\_3\_10} & \textit{inst\_4\_10} & \textit{inst\_5\_10} & \textit{inst\_6\_10} & \textit{inst\_7\_10} & \textit{inst\_8\_10} & \textit{inst\_9\_10} & \textit{inst\_10\_10} \\ \hline
		$W$ & \multicolumn{10}{c}{\bf Profit reduction ($\%$)}\\
		1 & 0\% &	$\infty$ &	59\%	& 4\%	& 68\%	& 15\%	&  $\infty$ &	26\%	& 17\%	& $\infty$\\
		2 & 0\% &	50\%	& 20\%	 & 4\%	& 68\%	& 15\%	& 55\%	& 26\%	& 17\%	& 36\%\\
		3 & 0\% &	50\%	& 20\%	& 4\%	& 68\% &	15\%	& 55\%	& 26\%	& 17\%	& 36\%\\
		4 & 0\%	& 0\%	& 0\% & 	4\%	& 68\%	& 36\%	& 0\%	& 26\%	& 17\%	& 0\%\\
		\hline
		$W$ & \multicolumn{10}{c}{\bf Computational time reduction ($\%$)}\\
		1 & 99\%	 & 87\%	& 100\%	& 99\%	& 99\%	& 100\%	 & 100\%	 & 100\%	 & 99\%	& 100\%\\
		2 & 95\% & 	56\%	& 97\%	 & 91\%	 & 97\%	& 99\%	& 97\%	& 95\%	& 94\%	& 99\%\\
		3 & 89\%	 & 0\% & 	91\%	& 79\%	& 86\%	& 94\%	& 89\%	& 86\%	& 84\%	& 97\%\\
		4 & 78\% & 	$-82\%$ & 	75\% & 	41\%	& 28\%	& 37\%	& 86\%	& 42\%	& $-151\%$ & 	82\%\\
		\bottomrule
\end{tabular}}
\caption{Detailed results on the performance of the rolling horizon approach, in terms of reduction of the profit and of the CPU time when compared to the $RP$ problem. The results refer to the instances with 10 bins.} \label{tab_RH_10bins}
\end{table}

\begin{table}[h!]
\centering
\resizebox{0.9\textwidth}{!}{
	\begin{tabular}{lllllllllll}\toprule
		& \textit{inst\_1\_11} & \textit{inst\_2\_11} & \textit{inst\_3\_11} & \textit{inst\_4\_11} & \textit{inst\_5\_11} & \textit{inst\_6\_11} & \textit{inst\_7\_11} & \textit{inst\_8\_11} & \textit{inst\_9\_11} & \textit{inst\_10\_11} \\ \hline
		$W$ & \multicolumn{10}{c}{\bf Profit reduction ($\%$)}\\
		1 & 49\% & 	93\%	 & 28\%	& 46\%	& 11\%	& 26\%	& 54\%	& 16\%	& 48\%	& 10\%\\
		2 & 49\%	 & 34\% & 	28\%	& 46\%	& 11\%	& 26\%	& 15\%	& 16\%	& 48\%	& 10\%\\
		3 & 49\%	& 34\%	& 28\%	& 46\%	& 11\%	& 26\%	& 15\%	& 16\%	& 48\%	& 10\%\\
		4 & 0\%	& 0\%	& 0\%	& 54\%	& 0\%	& 0\%	& 0\%	& 16\%	& 48\%	& 10\%\\
		\hline
		$W$ & \multicolumn{10}{c}{\bf Computational time reduction ($\%$)}\\
		1 & 100\%	 & 98\%	& 98\%	& 84\%	& 100\%	& 100\%	& 100\%	& 98\%	& 99\%	& 100\%\\
		2 & 99\%	 & 85\%	& 84\%	& 30\%	& 96\%	& 99\%	& 96\%	& 91\%	& 97\%	& 99\%\\
		3 & 95\%	& $-490\%$	& 61\% &	$-23$\%	 & 91\%	& 98\%	& 69\%	& 55\%	& 93\%	& 96\%\\
		4 & 82\%	 & $-141\%$ & 	$-36\%$	& $-161\%$ & 	80\%	& 93\%	& 66\%	& 10\%	& 10\%	& 90\%\\
		\bottomrule
\end{tabular}}
\caption{Detailed results on the performance of the rolling horizon approach, in terms of reduction of the profit and of the CPU time when compared to the $RP$ problem. The results refer to the instances with 11 bins.} \label{tab_RH_11bins}
\end{table}

\end{document}